\documentclass[12pt, oneside,reqno]{amsart}
\usepackage{Everything}
\usepackage[foot]{amsaddr}
\def\heit{{\mathsf{ht}}}

\renewcommand{\X}{\bar{X}}
\renewcommand{\Y}{\bar{Y}}
\renewcommand{\A}{\bar{A}}
\renewcommand{\B}{\bar{B}}

\let\oldtocsection=\tocsection

\let\oldtocsubsection=\tocsubsection

\let\oldtocsubsubsection=\tocsubsubsection

\renewcommand{\tocsection}[2]{\hspace{0em}\oldtocsection{#1}{#2}}
\renewcommand{\tocsubsection}[2]{\hspace{1em}\oldtocsubsection{#1}{#2}}
\renewcommand{\tocsubsubsection}[2]{\hspace{2em}\oldtocsubsubsection{#1}{#2}}

\makeatletter
\renewcommand\subsubsection{\@secnumfont}{\bfseries}
\renewcommand\subsubsection{\@startsection{subsubsection}{3}
  \z@{.5\linespacing\@plus.7\linespacing}{-.5em}%
  {\normalfont\bfseries}}

\title{Cartagena logic}

\author[F1.Kivimäki]{Siiri Kivimäki$^1$}
\address{$^1$Matematiikan ja tilastotieteen osasto, Helsingin yliopisto and Institut de Mathématiques de Jussieu - Paris Rive Gauche, Université Paris Cité.}
\author[F2.Väänänen]{Jouko Väänänen$^2$}
\address{$^2$Department of Mathematics and Statistics, University of Helsinki and ILLC, University of Amsterdam}
\author[F3.Villaveces]{Andrés Villaveces$^3$}
\address{$^3$Departamento de Matemáticas, Universidad Nacional de Colombia - Bogotá}

\begin{document}

\maketitle

\begin{abstract}
    We introduce a new kind of infinitary logic that we call \textit{Boolean expansion of $\mathcal{L}_{\kappa\kappa}$}. This logic involves a new kind of variable, that we call \textit{generalised Boolean variable}. These variables range over the powerset of a cardinal number in a way reminiscent of random variables. From this Boolean expansion, we extract a traditional infinitary logic, called Cartagena logic. We prove several model-theoretic properties of Cartagena logic, and give multiple examples of its expressive power. The main result is that Cartagena logic is a good syntactically defined approximation to Shelah's infinitary logic $\mathcal{L}^1_\kappa$ (from \cite{shelah2012nice}). The latter is not known to have a generative syntax, while Cartagena logic does have a very clear one.
\end{abstract}

\tableofcontents

\section*{Introduction}

\subsection*{A syntactic approach to Shelah's logic}\spa 

\vv

The combination of Löwenheim-Skolem Theorem and Compactness Theorem limits the expressive power of a logic to that of first order logic. This maximality principle is the famous Lindström Theorem for first order logic (Lindström \cite{lindstrom1969}). It reveals that first order logic is at an optimal point of balance: by adding expressive power to it one necessarily loses model-theoretic properties. Soon after Lindström's result, the question was raised, whether there are other logics at a similar point of equilibrium. More precisely: are there strict strengthenings of first order logic satisfying a Lindström-type characterization? Despite the naturality of this question, it remained unanswered, until recently.

In 2012, Shelah \cite{shelah2012nice} offered a solution to this problem, in the form of a logic he calls $\LL^1_\kappa$ (where $\kappa$ is an uncountable cardinal with $\kappa=\beth_\kappa$). The logic $\LL^1_\kappa$ is an infinitary logic strictly between the logics $\LL_{\kappa\omega}$ and $\LL_{\kappa\kappa}$. It has a Lindström-type characterization in terms of a property called \textit{strong well ordering number $\kappa$}, which is a mode-theoretic property combining weak forms of Compactness and Löwenheim-Skolem type of properties (see Definition \ref{sudwodef} and Theorem \ref{firstcha}).

In all known cases, a proof of a Lindström-type characterization simultaneously gives a proof of interpolation\footnote{This said, it is still open even for first order logic whether one can replace compactness by interpolation in Lindström theorem. For definition of interpolation, see Remark \ref{inter}.}. This is the case for $\LL^1_\kappa$ too, which balances the classical result of Malitz \cite{malitz1971infinitary} of interpolation for $\LL_{\kappa\omega}$, happening  in $\LL_{\kappa\kappa}$.

In addition to interpolation and strong well ordering number $\kappa$, each complete $\LL^1_\kappa$-theory admits some sort of special models (Shelah \cite{shelah2021isomorphic}), and for $\kappa$ strongly compact, its elementary equivalence has an algebraic characterization in the spirit of classical Keisler-Shelah Theorem for first order logic. Namely: two models are $\LL^1_\kappa$-elementary equivalent if and only if they have isomorphic iterated ultrapowers along a countable sequence of $\kappa$-complete ultrafilters (Shelah \cite{shelah2021isomorphic}). The logic $\LL^1_\kappa$ has thus proved quite promising when it comes to its model theory. 

There is, however, one aspect where $\LL^1_\kappa$ seems to be rather weak: the syntax. The logic $\LL^1_\kappa$ is derived from a game, in the sense that a sentence is, by definition, a class of structures, closed under a certain Ehrenfeucht-Fraïssé type of game. This results in the absence of a \textit{generative syntax}, i.e. a syntax defined in such a way that the set of all formulas can be obtained by closing the set of atomic formulas under negation, conjunction, quantifiers, and possibly other logical operations.

The lack of generative syntax complicates, on one hand, applying the great model-theoretic properties of $\LL^1_\kappa$ in real life situations, and on the other hand, further study of it and logics in its neightborhood. It is clear that many classical methods, such as the method of Skolem functions or simply induction on the complexity of formula, highly requires generative syntax. We believe that much more can be said about $\LL^1_\kappa$ and even about more general questions regarding logics derived from games, by providing a generative syntax for $\LL^1_\kappa$.

In the present paper, we address the general question of deriving a syntax from a game, and the more localized question of finding a syntax for $\LL^1_\kappa$. Partial answers are provided to both questions.

Our approach is the following: 
\begin{itemize}
    \item We first define an expansion $\LL^{\Bool}_{\kappa\kappa}$ of $\LL_{\kappa\kappa}$ that we call the \textit{Boolean expansion} of $\LL_{\kappa\kappa}$, which has exactly the same sentences as $\LL_{\kappa\kappa}$, but in which a new kind of variable is allowed in formulas. We believe that this expansion is useful when approaching the general question of deriving a syntax from a game. The reason is that these new variables allow, in a  sense, quantifying over moves in a game (see more precise discussion in the beginning of section \ref{upwards}).
    \item We then define a game, Cartagena game, which is a simplified variant of the game of $\LL^1_\kappa$. With the help of $\LL^{\Bool}_{\kappa\kappa}$, we are able to build a syntax that completely corresponds to Cartagena game.
    \item Finally we study the model-theory of Cartagena logic. It turns out that some model-theoretic properties can be proved in a stronger form than in $\LL^1_\kappa$. This, as expected, results in a slightly weaker expressive power. 
    \item However, we are able to show that Cartagena logic is not too much weaker than $\LL^1_\kappa$, in a way that can be made precise by means of the $\Delta$-operator. Also, Cartagena logic does not have interpolation nor Lindström Theorem, contrary to $\LL^1_\kappa$.
\end{itemize}
Our main result, thus, is the existence of a good approximation to logic $\LL^1_\kappa$, with a simple generative syntax and rich model theory.

\subsection*{Organization of the paper and notation}

This is a paper in abstract model theory. For introduction to the subject, see \cite{Barwise1985-BARML-8}. For context, see \cite{shelah2012nice}, \cite{shelah2021isomorphic} and \cite{dvzamonja2021chain}.

We deal with infinitary logics, in particular those whose expressive power fall between the logic $\LL_{\kappa\omega}$, allowing conjunctions and disjunctions of length $<\kappa$, and the logic $\LL_{\kappa\kappa}$, allowing in addition quantification over tuples of length $<\kappa$. We assume familiarity with these logics.

The first section is devoted to the definition of Cartagena logic: Cartagena game, Cartagena syntax, and finally the proof of game-syntax correspondence. In the second section we prove three major model theoretic properties of Cartagena logic: closure under countable elementary chains (Union Lemma \ref{union}), Löwenheim-Skolem-Tarski Theorem (\ref{LST}) and a strong form of undefinability of well order (\ref{sudwo}).
The last section describes the expressive power of Cartagena logic, first through explicit examples, and then by comparing it with $\LL^1_\kappa$.

Throughout, $\kappa$ denotes an uncountable cardinal. In sections \ref{mtocl} and \ref{sectexp} it is assumed in addition that $\kappa$ satisfies $\kappa=\beth_\kappa$. In case $\kappa$ is a singular cardinal, the logics $\LL_{\kappa\omega}$ and $\LL_{\kappa\kappa}$ are defined by
\begin{align*}
    & \LL_{\kappa\omega}:=\bigcup_{\lambda<\kappa}\LL_{\lambda^+\omega},\\
    &\LL_{\kappa\kappa}:=\bigcup_{\lambda<\kappa}\LL_{\lambda^+\lambda^+}.
\end{align*}

Our notation is mostly standard, with the following exception:
For variables $x_i$ or elements $a_i$ and an index set $A$, we denote \[
\x_A:=(x_i)_{i\in A}\quad\text{and}\quad \a_A:=(a_i)_{i\in A}.
\]

We follow the convention that
\begin{align*}
    &\bigwedge_{\emptyset}\quad:=\quad\top,\\
    &\bigvee_{\emptyset}\quad:=\quad\bot.
\end{align*}
A \textit{signature} is a set of relation, function and constant symbols of finite arity. In our paper, for simplicity, we assume that a signature is always relational and of size $<\kappa$. We tacitly fix such a signature and won't explicitly mention it.

\subsection*{Acknowledgements} We wish to thank Will Boney, Xavier Caicedo, Mirna Džamonja and Boban Veličkovi\'c for their comments on earlier versions of this work. This project has received funding from the European Research Council (ERC) under the
European Union’s Horizon 2020 research and innovation programme (grant agreement No
101020762). The second author was supported by the Academy of Finland grant number 322795.

\section{Boolean expansion $\LL^{\Bool}_{\kappa\kappa}$}

In this section, an expansion $\LL^{\Bool}_{\kappa\kappa}$ of the logic $\LL_{\kappa\kappa}$ is built, for each infinite cardinal $\kappa$. This expansion $\LL^{\Bool}_{\kappa\kappa}$ allows formulas with two types of variables: the ordinary ones and ``Boolean'' ones, that range over the powerset of a cardinal.

\subsection{Formulas with Boolean variables}\spa 

\vv

\begin{de}
    A \textit{(generalized) Boolean variable} is a variable $X$ which ranges over subsets of some cardinal $\theta$, i.e. over the Boolean algebra $\Pow(\theta)$. 
\end{de}

Generalized Boolean variables are not quite the same as second order variables: they do not range over subsets of an (intended) domain (in the interpretation); rather, they range over subsets of the cardinal $\theta$.

We make the idea more precise in what follows. First we will fix sets of standard variables $\Var_\kappa$ and Boolean variables $\BVar_\kappa$, for each cardinal $\kappa$, that are intended to be used in formulas of the Boolean extendion $\LL_{\kappa\kappa}^{\Bool}$:

\begin{notation}\label{variables} \spa 
    \begin{enumerate}
        \item For each cardinal $\lambda$, we fix a set $\mathsf{Var}_\lambda$ of size $\lambda$ of \textit{standard variables} intended to range over elements of model.
        \item For each cardinal $\lambda$ and each cardinal $\theta<\lambda$, we fix a set $\mathsf{BVar}_\lambda^\theta$ of size $\lambda$ of \textit{Boolean variables}, intended to range over $\Pow(\theta)$.
        \item For each Boolean variable $X\in\BVar_\lambda^\theta$, we fix a $\theta$-tuple of pairwise distinct standard variables
        \[
        \x_\theta^X\subseteq\Var_\lambda.
        \]
        \item We write
        \[
        \BVar_\lambda:=\bigcup_{\theta<\lambda}\BVar_\lambda^\theta,.
        \]
    \end{enumerate}
    Whenever ${\lambda'}\leq\lambda$, we assume 
    \begin{align*}
        &\Var_{\lambda'}\subseteq\Var_\lambda,\\
        &\BVar_{\lambda'}^\theta\subseteq\BVar_\lambda^\theta.
    \end{align*}
    We also assume that the sets of Boolean variables $\BVar_\lambda^\theta$ as well as the corresponding sets of standard variables $\{x_i^X:i\in\theta\}$ are pairwise disjoint.

    When $\lambda$ is clear from context, we abbreviate by writing $\Var$ for $\Var_\lambda$, $\BVar^\theta$ for $\BVar_\lambda^\theta$ and $\BVar$ for $\BVar_\lambda$.
    \end{notation}

A Boolean variable will appear in the place of an index set in conjunctions and disjuntions. In order to emphasize a distinction between a ``real'' disjunction and a disjunction with a variable in the place of an index set, we will make use of the symbols 
\begin{align*}
    &\bigdvee,\\
    &\bigdwedge.
\end{align*}
Then, when a real set is substituted for the variable, these symbols will be interpreted as disjunction and conjunction, respectively.

We now expand the set of $\LL_{\kappa\kappa}$-formulas to contain formulas with Boolean variables. Some examples and explanations are given right after the definition. We tacitly fix a relational signature. An \textit{atomic formula} is then a formula of the form
\[
    R(\x),
\]
where $R$ is either the symbol $=$, or an $n$-ary relation symbol of the signature, and $\x$ is a finite tuple of standard variables.

\begin{de}[Boolean extension $\LL^{\Bool}_{\kappa\kappa}$]\label{boolext} Let $\lambda$ be an infinite regular cardinal.
    The \textit{$\LL_{\lambda\lambda}^{\mathsf{Bool}}$-formulas} are defined as follows.
    \begin{enumerate}
        \item Every atomic formula with variables in $\Var_\lambda$ is an $\LL_{\lambda\lambda}^{\mathsf{Bool}}$-formula.
        \item If $\phi$ is an $\LL_{\lambda\lambda}^{\mathsf{Bool}}$-formula then so is
        \begin{align*}
            &\neg\phi.
        \end{align*}
        \item If $\Phi$ is a set of $\LL_{\lambda\lambda}^{\mathsf{Bool}}$-formulas of size $<\lambda$, then 
        \begin{align*}
            &\bigvee\Phi,\\
            & \bigwedge\Phi,
        \end{align*}
        are $\LL_{\lambda\lambda}^{\mathsf{Bool}}$-formulas.
        \item If $\phi$ is an $\LL_{\lambda\lambda}^{\mathsf{Bool}}$-formula and $\x\subseteq\mathsf{Var}_\lambda$ has length $<\lambda$, then \begin{align*}
            &\exists\x\phi,\\
            & \forall\x\phi,
        \end{align*} 
        are $\LL^{\Bool}_{\lambda\lambda}$-formulas.
        \item\label{batomic} If $\{\phi_u:u\subseteq\theta\}$ are $\LL_{\lambda\lambda}^{\mathsf{Bool}}$-formulas, then so are the following:
        \begin{align*}
            &\bigdvee_{u\in p(X)}\phi_u,\\
            &\bigdwedge_{u\in p(X)}\phi_u,
        \end{align*}
        provided that $\theta<\lambda$, $X\in\mathsf{BVar}^\theta_\lambda$ and $p:\Pow(\theta)\to\Pow(\Pow(\theta))$ is a function.
        \item\label{unarq} If $\phi$ is an $\LL^{\Bool}_{\lambda\lambda}$-formula, then so are the following:
        \begin{align*}
            &\bigdvee_{X\in\WW}\phi,\\
            &\bigdwedge_{X\in\WW}\phi,
        \end{align*}
        provided that $\theta<\lambda$ and $\WW\subseteq\Pow(\theta)$.
        \item\label{longq} If $\phi$ is an $\LL^{\Bool}_{\lambda\lambda}$-formula, $\X=(X_i)_i\subseteq\BVar_\lambda$ is a tuple of Boolean variables of length $<\lambda$ with $X_i\in\BVar_\lambda^{\theta_i}$ and $\WW_i\subseteq\Pow(\theta_i)$, then
        \begin{align*}
            &\bigdvee_{\X\in\prod_i\WW_i}\phi,\\
            &\bigdwedge_{\X\in\prod_i\WW_i}\phi
        \end{align*}
        are $\LL^{\Bool}_{\lambda\lambda}$-formulas.
    \end{enumerate}
    For a singular cardinal $\kappa$, we let 
    \[
    \LL_{\kappa\kappa}^{\Bool}:=\bigcup_{\lambda<\kappa}\LL^{\Bool}_{\lambda^+\lambda^+}.
    \]
\end{de}

A Boolean variable can occur only as an (uninterpreted) index set of a conjunction or a disjunction. The operation
\begin{align*}
    &\bigdvee_{u\in p(X)}
\end{align*}
from clause \ref{batomic} is a logical operation that introduces $X$ as a variable in the formula. It is not a ``real'' disjunction. Later, we will see how to substitute a set $A$ for $X$, resulting in a ``real'' disjunction
\[
\bigvee_{u\in p(A)}.
\]
The definition for substitution is given below, Definition \ref{substitution}. Similarly for the conjunction $\bigdwedge_{u\in p(X)}$.

\begin{ex}\label{ex}
    For example, the formula
\begin{align*}
    & \bigdvee_{i\in X}P(x_i^X)
\end{align*}
is a formula with one free Boolean variable $X$ and one tuple of free standard variables $\x_\theta^X$, in signature containing a unary predicate symbol $P$. When a set $A\subseteq\theta$ is substituted for $X$, the result is
\begin{align*}
    & \bigvee_{i\in A}P(x_i^X).
\end{align*}
More examples are given in \ref{examples}.
\end{ex}

Clause \ref{unarq} corresponds to a ``unary'' quantification of a Boolean variable, and clause \ref{longq} corresponds to quantifying over long tuples of Boolean quantifiers.

\begin{rmk}
    The construction of $\LL^{\Bool}_{\kappa\kappa}$-formulas can be formalized within set theory in a standard way, as on page 6 in \cite{Keisler1971}, or as Definition 1.3 on page 81 in \cite{barwise2017admissible}.
\end{rmk}

\subsection{Free Boolean variables and substitution}\spa 

\vv

Before proceeding further, we discuss free Boolean variables. The definition is natural but subtle. For the sake of notational simplicity, we use the notation $\x$ both for the tuple and for its range $\{x_i\}_i$, and similarly for tuples $\X$ of Boolean variables.

\begin{de}
    Let $\kappa$ be an infinite cardinal and let $\phi$ is an $\LL^{\Bool}_{\kappa\kappa}$-formula. With the notation from Definition \ref{boolext}, we define:
    \begin{enumerate}
        \item The set of \textit{free (standard) variables of $\phi$,} $\FreeVar(\phi)\subseteq\Var_\kappa$ is defined as follows:
        \begin{align*}
            &\FreeVar\left(\phi\right):=\{x\in\Var_\kappa:x\text{ occurs in }\phi\},\quad \text{for }\phi\text{ atomic,}\\
            & \FreeVar\left(\neg\phi\right):=\FreeVar\left(\phi\right),\\
            & \FreeVar\left(\bigvee\Phi\right):=\bigcup_{\phi\in\Phi}\FreeVar\left(\phi\right),\\
            &\FreeVar\left(\exists\x\phi\right):=\FreeVar(\phi)-\x,\\
            &\FreeVar\left(\bigdvee_{u\in p(X)}\phi_u\right):=\bigcup_{u\subseteq\theta}\FreeVar\left(\phi_u\right),\\
            &\FreeVar\left(\bigdvee_{X\in\WW}\phi\right):=\FreeVar(\phi),\\
            &\FreeVar\left(\bigdvee_{\X\in\prod_i\WW_i}\phi\right):=\FreeVar(\phi).
        \end{align*}
        Similarly for the dual clauses given by $\forall$ and $\bigwedge$.
        \item The set of \textit{free Boolean variables of $\phi$}, $\FreeBVar(\phi)\subseteq\BVar_\kappa$ is defined as follows:
        \begin{align*}
            &\FreeBVar\left(\phi\right):=\emptyset\quad \text{for }\phi\text{ atomic,}\\
            & \FreeBVar\left(\neg\phi\right):=\FreeBVar\left(\phi\right),\\
            & \FreeBVar\left(\bigvee\Phi\right):=\bigcup_{\phi\in\Phi}\FreeBVar\left(\phi\right)\\
            &\FreeBVar\left(\exists\x\phi\right):=\FreeBVar(\phi)\\
            &\FreeBVar\left(\bigdvee_{u\in p(X)}\phi_u\right):=\bigcup_{B\subseteq\theta}\FreeBVar\left(\phi_u\right)\cup\{X\},\\
            &\FreeBVar\left(\bigdvee_{X\in\WW}\phi\right):=\FreeBVar(\phi)-\{X\},\\
            &\FreeBVar\left(\bigdvee_{\X\in\prod_i\WW_i}\phi\right):=\FreeBVar(\phi)-\X.
        \end{align*}
        Similarly for the dual clauses given by $\forall$ and $\bigwedge$.
    \end{enumerate}
\end{de}

\begin{de} Let $\kappa$ be an infinite cardinal and let $\phi$ be an $\LL^{\Bool}_{\kappa\kappa}$-formula.
    Let $\x=(x_i)_i$ be a tuple of variables from $\Var_\kappa$ and let $\bar{X}=(X_i)_i$ be a tuple of variables from $\BVar_\kappa$. We use the notation
    \[
    \phi(\x,\bar{X})
    \]
    to indicate that
    \begin{enumerate}
        \item $\FreeVar(\phi)\subseteq\{x_i\}_i$, and
        \item $\FreeBVar(\phi)\subseteq\{X_i\}_i$.
    \end{enumerate}
\end{de}

Substituting a tuple of elements $\a\subseteq\MM^{\x}$ for a tuple of variables $\x$ in a formula $\phi(\x,\X)$ is done in a standard way, as in the case of the logic $\LL_{\kappa\kappa}$. The next goal is to substitute a set $A$ for a Boolean variable $X$, and finally see how every $\LL^{\Bool}_{\kappa\kappa}$-formula $\phi(\x,\X)$ can be ``collapsed'' to an $\LL_{\kappa\kappa}$-formula via a substitution $\phi(\x,\X)\mapsto\phi(\x,\X)[\A/\X]$. See Definition \ref{substitution} and Proposition \ref{collapsing}.

\begin{rmk}\label{once}
    Up to a change of symbols, we may assume that every Boolean variable $X$ that occurs free in an $\LL^{\Bool}_{\kappa\kappa}$-formula, occurs exactly once. Furthermore, we may assume that no variable that occurs free occurs also bounded in one formula. From now onwards, we tacitly assume that this always happens.
\end{rmk}

\begin{de}
    Let $\X=(X_i)_i$ be a tuple of $\BVar_\kappa$, with $X_i\in\BVar^{\theta_i}_\kappa$. A \textit{valuation for $\X$} is a tuple $\A=(A_i)_i$ of same domain such that 
    \[
    A_i\subseteq\theta_i
    \]
    for each $i$.
\end{de}

\begin{rmk}
    In particular: any element of a set $\WW\subseteq\Pow(\theta)$ is a valuation for any $X\in\BVar^{\theta}_\kappa$. Similarly for longer tuples $\X$ and $\A\in\prod_i\WW_i$.
\end{rmk}

In the next definition we write $X\in\bar{X}$ instead of $X\in\{X_i\}_i$, when no confusion arises.

\begin{de}\label{substitution} Let $\phi(\x,\X)$ be an $\LL^{\Bool}_{\kappa\kappa}$-formula, let $\X=(X_i)_i$ be a tuple of Boolean variables from $\BVar_\kappa$, let $\X'$ be a subtuple of $\X$, and let $\A'$ be a valuation for $\X'$.
We define the formula
\[
    \phi(\x,\X)[\A'/\X']
\]
inductively, as follows:
\begin{align*}
        \phi(\x,\X)[\A'/\X']&:=\phi(\x,\X),\quad\text{for atomic }\phi(\x),\\
        \neg\phi(\x,\X)[\A'/\X']&:=\neg(\phi(\x,\X)[\A'/\X']),\\
        \bigvee\Phi(\x,\X)[\A'/\X']&:=\bigvee_{\phi\in\Phi}\phi(\x,\X)[\A'/\X'],\\
        (\exists\y\phi)(\y,\x,\X)[\A'/\X' ]&:=\exists\y\phi(\y,\x,\X)[\A'/\X'],\\
        \left(\bigdvee_{u\in p(X_i)}\phi_u(\x,\X)\right)[\A'/\X' ]&:=\begin{cases}
            \bigvee_{u\in p(A_i)}\phi_u(\x,\X)[\A'-\{A_i\}/\X'-\{X_i\}], \quad &\text{if }X_i\in \X',\\
            \bigdvee_{u\in p(X_i)}\phi_u(\x,\X) [\A'/\X'], &\text{if }X_i\notin\X',
        \end{cases}\\
        \left(\bigdvee_{Y\in \WW}\phi(\x,Y,\X)\right)[\A'/\X' ]&:=\bigvee_{B\in\WW}\phi(\x,Y,\X)[B/Y][\A'/\X'],\\
      \left(\bigdvee_{\Y\in\prod_j\WW_j}\phi(\x,\Y,\X)\right)[\A'/\X' ]&:=\bigvee_{\B\in\prod_j\WW_j}\phi(\x,\Y,\X)[\B/\Y][\A'/\X'].
\end{align*}
Similarly for the dual clauses for conjunctions and universal quantifier.
\end{de}

\begin{rmk}
    The crucial clauses are the one with the operation $\bigdvee_{u\in p(X)}$ and the one with $\bigdvee_{Y\in\WW}$. In the simplest cases these clauses become:
    \begin{align*}
        &\left(\bigdvee_{u\in p(X)}\phi_u\right)[A/X]=\bigvee_{u\in p(A)}\phi_u,\\
        &\left(\bigdvee_{Y\in\WW}\phi\right)=\bigvee_{B\in\WW}\phi[B/Y].
    \end{align*}
\end{rmk}

\begin{notation}
    For an $\LL^{\Bool}_{\kappa\kappa}$-formula $\phi(\x,\X)$ and a valuation $\A$ for $\X$, we define
    \[
    \phi(\x,\A):=\phi(\x,\X)[\A/\X].
    \]
\end{notation}

\subsection{Projection to $\LL_{\kappa\kappa}$}\spa 

\vv

For reasons of presentation, we introduce the following auxiliary concept: A \textit{global valuation} in $\LL^{\Bool}_{\kappa\kappa}$ is a tuple 
\[
\A=(A_X)_{X\in\BVar_\kappa}
\]
such that $A_X\subseteq\theta$ whenever $X\in \BVar^{\theta}_\kappa$. Thus, a global valuation is a simultaneous valuation for every $X\in\BVar_\kappa$. The following proposition immediately follows from definitions:

\begin{prop}\label{collapsing}
    Every global valuation $\A$ uniquely determines a mapping
    \[
\begin{tikzcd}
        \LL^{\Bool}_{\kappa\kappa}\text{-formulas}\arrow{r}{\pi_{\bar{A}}}&\LL_{\kappa\kappa}\text{-formulas}
\end{tikzcd}
\]
via the substitution
\[
\begin{tikzcd}
        \phi(\x,\X)\arrow[|->]{r} & \phi(\x,\A).
\end{tikzcd}
\]
\end{prop}

With valuations we get semantics:

\begin{de}\label{semantics}
    Let $\phi(\x,\X)$ be an $\LL^{\Bool}_{\kappa\kappa}$-formula. For a structure $\MM$, a tuple $\a\in\MM^{\x}$ and a valuation $\A$ of $\X$, we denote
    \[
    \MM\models \phi(\a,\A)
    \]
    if $\MM\models\phi(\x,\A)[\a/\x]$, where $\phi(\x,\A)$ is understood as an $\LL_{\kappa\kappa}$-formula via Proposition \ref{collapsing}, and $[\a/\x]$ denotes the ordinary substitution of $\a$ for $\x$, defined for the logic $\LL_{\kappa\kappa}$.
\end{de}

In conclusion:

\begin{rmk}\label{collapse} There are two different ways of transforming an $\LL^{\Bool}_{\kappa\kappa}$-formula into an $\LL_{\kappa\kappa}$-formula. 
\begin{enumerate}
    \item \textit{Valuation:} Given an $\LL^{\Bool}_{\kappa\kappa}$-formula $\phi(\x,\X)$ and a valuation $\A$ for $\X$, the formula
    \[
    \phi(\x,\A)
    \]
    is equivalent to an $\LL_{\kappa\kappa}$-formula.
    \item \textit{Quantification:} Given an $\LL^{\Bool}_{\kappa\kappa}$-formula $\phi(\x,\X)$, where $X_i\in\BVar^{\theta_i}_\kappa$ and sets $\WW_i\subseteq\Pow(\theta_i)$, the formulas
    \begin{align*}
        & \bigdvee_{\X\in\prod_i\WW_i}\phi(\x,\X),\\
        & \bigdwedge_{\X\in\prod_i\WW_i}\phi(\x,\X)
    \end{align*}
    are equivalent to $\LL_{\kappa\kappa}$-formulas. 
\end{enumerate}
\end{rmk}

An $\LL^{\Bool}_{\kappa\kappa}$-\textit{sentence} is an $\LL^{\Bool}_{\kappa\kappa}$-formula $\phi$ with no free variables:
\[
    \FreeVar(\phi)=\FreeBVar(\phi)=\emptyset.
\]

\begin{rmk}
    Every $\LL^{\Bool}_{\kappa\kappa}$-formula with no free Boolean variables is (equivalent to) an $\LL_{\kappa\kappa}$-formula. In particular, every $\LL^{\Bool}_{\kappa\kappa}$-sentence is (equivalent to) an $\LL_{\kappa\kappa}$-sentence. This follows from the previous remark. Thus, $\LL^{\Bool}_{\kappa\kappa}$ and $\LL_{\kappa\kappa}$ are the same logic, when it comes to (sentential) expressive power:
    \[
    \LL^{\Bool}_{\kappa\kappa}\equiv\LL_{\kappa\kappa}.
    \]
    Formulas with free variables distinguish $\LL^{\Bool}_{\kappa\kappa}$ from $\LL_{\kappa\kappa}$.
\end{rmk}

\section{Definition of Cartagena logic}   

Shelah's logic $\LL^1_\kappa$ is derived from an Ehrenfeucht-Fraïssé type of game, and Cartagena logic approximates $\LL^1_\kappa$. To define Cartagena logic, it is thus necessary to start from a game, a simplified version of Shelah's game. Indeed, the goal of this section is to first give the definition of Cartagena game (Definition \ref{game}), then give the definition of Cartagena syntax (Definition \ref{syntax}), and finally prove that they correspond to each other in the following sense: player $\2$ has a winning strategy in the Cartagena game on two structures if and only if the structures satisfy the same Cartagena sentences. This is Theorem \ref{gamesyntax}. 

\subsection{Cartagena game}\spa 

\vv

Cartagena game is a simplification of Shelah's game for $\LL^1_\kappa$. As is the latter, also Cartagena game is an Ehrenfeucht-Fraïssé type of game played by two players on two structures, where the job of player $\2$ is to approximate an isomorphism between the two models. A game like this gives rise to a logic, roughly, as follows: Firstly, it gives rise to a notion of elementary equivalence, via two structures claimed to be elementarily equivalent if player $\2$ has a winning strategy in the game on them. Secondly, sentences can be taken to be unions of equivalence classes, where a model satisfies a sentence, by definition, if it belongs to it as an element. See definition of $\LL^1_\kappa$ in section \ref{sectshlogic} for a concrete example. This approach is natural in the sense that a sentence, by definition, becomes the class of models it defines.

We will give the definition of Cartagena game, leaving out the derivation of the logic from it, as, in section \ref{sectsyntax}, an \textit{explicit, generative, syntactic} definition of Cartagena logic is provided.

We first describe the game informally. In the usual Ehrenfeucht-Fraïssé game for $\LL_{\kappa\kappa}$, there are two models and two players. Player $\1$ picks tuples of length $<\kappa$, and player $\2$ has to map them to the opposite model. Player $\1$ also has a clock - he has to play ordinals below $\kappa$ in a strictly descending order. This, of course, results the game to end after finitely many steps. The idea in Cartagena game (and Shelah's \textit{delayed} game) is similar. The clock is there, and player $\1$ plays tuples of size $<\kappa$, which player $\2$ has to map in the opposite model. However, she does not need to map the whole tuple. Instead, she partitions the set into countably many pieces, and has to be ready to map any one of the pieces.

We give now the definition in detail.

\vv 

\begin{de}\label{partition} Let $\theta$ be an ordinal and  $f:\theta\to\omega$  a function. We write
    \[
        \WW_f:=\{f^{-1}\{n\}:n\in\omega\}
    \]
and call this set the \textit{partition of $\theta$ given by $f$}.
\end{de}

\begin{de}[Cartagena game $\G^\beta_\lambda$]\label{game} Let $\lambda$ be a cardinal. Let $\MM$ and $\NN$ be structures in a same signature and let $\beta$ be an ordinal. We define the \textit{Cartagena game} of height $\beta$
\[
\G_\lambda^\beta(\MM,\NN).
\]
The \textit{states} of the game are pairs $(\alpha,\pi)$, where $\alpha\leq\beta$ is an ordinal and $\pi:\MM\to\NN$ is a partial isomorphism.

\textbf{Starting state:} The starting state is $(\beta,\emptyset)$.

\textbf{Further states:} At state $(\alpha,\pi)$: 
    \begin{enumerate}
        \item Player $\1$ picks an ordinal $\alpha'<\alpha$ and a tuple $\a_\theta\in\MM^\theta$ (or a tuple $\b_\theta\in\NN^\theta$), for some cardinal $\theta<\lambda$.
        \item Player $\2$ picks a function $f:\theta\to\omega$ and a tuple $\b_\theta\in\NN^\theta$
        (or $\a_\theta\in\MM^\theta$), such that for any $A\in \WW_f$, the map
        \[
        \pi_A:=\pi\cup\{(a_i,b_i):i\in A\}
        \]
        is a partial isomorphism.
        \item Player $\1$ picks $A\in\WW_f$.
    \end{enumerate}
    
    The next state is $(\alpha',\pi_A)$.
    
\noindent The player who first cannot move loses.

\end{de}

A strategy for player $\2$ in the game $\G_\lambda^\beta(\MM,\NN)$ is a function $\sigma$ that takes in a state $(\alpha,\pi)$ and a move $(\alpha',\a_\theta)$ of player $\1$ and gives a move $(f,\b_\theta)$ of player $\2$. We end the subsection by stating a lemma, which will be important in the proof of game-syntax correspondence (Theorem \ref{gamesyntax}) and which shows that player $\2$ is allowed to play a finer function while preserving her winning strategy. 

We say that a partition $\WW$ \textit{refines} another partition $\VV$ if for every $A\in \WW$ there is $B\in\VV$ such that $A\subseteq B$.

\begin{de}
For functions $f,g:\theta\to\omega$, we write
\[
    f\leq g\quad:\iff\quad \WW_f\text{ refines }\WW_g.
\]
\end{de}

The following lemma is an immediate consequence of the definitions:
\begin{lem}\label{strategy}
    Let $f,g:\theta\to\omega$ be functions such that $f\leq g$. Suppose that there is a winning strategy $\sigma$ for player $\2$ in a game $\G_\lambda^\beta(\MM,\NN)$ such that
    \[
    \sigma((\alpha,\pi),(\alpha',\a))=(g,\b),
    \]
    for a state $(\alpha,\pi)$ and a move $(\alpha',\a)$ of player $\1$.
    Then there is a winning strategy $\sigma'$ for player $\2$ which agrees with $\sigma$ up to\footnote{I.e. for all states $s=(\beta,\pi')$ where $\beta>\alpha$ and for all moves $m$ of player $\1$, $\sigma(s,m)=\sigma'(s,m)$.} state $(\alpha,\pi)$ but
    \[
    \sigma'((\alpha,\pi),(\alpha',\a))=(f,\b).
    \]
\end{lem}

\subsection{Syntax of Cartagena logic}\label{sectsyntax}\spa 

\vv

This section is devoted to defining \textit{Cartagena syntax}, a syntax that corresponds to Cartagena game.

The definition involves two distinct steps:
\begin{itemize}
    \item First, we extract two important classes of formulas: \textit{upwards correct formulas} and \textit{downwards correct formulas}. These are classes of formulas that have a kind of monotonicity property with respect to substitution for a Boolean variable.
    \item Then, Cartagena syntax is defined with the help of these subclasses. The upwards and downwards correct formulas are the formulas that can be closed into a sentence with a ``Cartagena quantifier''.
\end{itemize} 
The first step is the most complicated part. Our earlier (simpler) attempts to capture the syntax of Shelah's logic made well order definable (in subtle ways) and were therefore wrong, as well order cannot be definable in any syntax corresponding to the Cartagena game (a fact that can quite easily be inferred from the definition of the game). Our construction will culminate in the definition of Cartagena syntax, Definition \ref{syntax}.

We start by upwards and downwards correctness.

\subsubsection{Upwards and downwards correctness}\label{upwards}\spa 

\vv 

The notions of \textit{upwards} or \textit{downwards correctness} are crucial when extracting the set of Cartagena sentences from the set of $\LL_{\kappa\kappa}$-sentences. By scrutinizing the Cartagena game, it is possible to prove that Cartagena logic lies strictly between $\LL_{\kappa\omega}$ and $\LL_{\kappa\kappa}$. In fact, it is much closer to $\LL_{\kappa\omega}$ than to $\LL_{\kappa\kappa}$. One of its features is that it is too weak to define the concept of  well ordering. Namely, the $\LL_{\kappa\kappa}$-sentence
\[
\neg\exists\x_\omega\bigwedge_{n\in\omega}x_{n+1}<x_n,
\]
which defines the class of well founded orders, is not (equivalent to) a Cartagena sentence. This can be seen by looking at the game (see Section \ref{sectsudwo} about undefinability of well order). Since the sentence defining well foundedness has a form that is almost the simplest possible among sentences in $\LL_{\kappa\kappa}$ but not in $\LL_{\kappa\omega}$, we were pushed to develop a kind of ``delayed conjunction" bounded by a kind of ``randomized long existential quantifier''. These quantifiers come in two dual forms, one of which looks like this:
\[
\exists\x_\theta\bigvee_{f:\theta\to\omega}\bigdwedge_{X\in\WW_f},
\]
where $\WW_f=\{f^{-1}\{n\}:n\in\omega\}$.
Furthermore, in Cartagena logic this quantifier can be added in front of \textit{only some kind of formulas, not in front of all of them}. This idea is made clear with the help of our concept of ``upwards or downwards correct'' formula. The set of upwards or downwards closed formulas will be the formulas that can be closed into a sentences by the above kind of ``Cartagena quantifier''. 

What pushed us to come up with Boolean variables was the problem of incorporating into the syntax the act of partitioning a set into countably many pieces, by player $\2$ in Cartagena game. In a sense, while a Boolean variable ranges over subsets of some $\theta$, in effect, it ranges over subtuples of some long tuple $\a_\theta$ played by player $\1$. Then, when a Boolean variable $X\in\BVar^\theta_\kappa$ takes a value $A\subseteq\theta$, this, in some sense, corresponds to fixing the subtuple
\[
\a_A\subseteq\a_\theta.
\]
By integrating Boolean variables into the syntax, we allow formulas to be undecided about player $\1$'s choice concerning which piece of partition should the game be continued with.

However, integrating Boolean variables turned out not to be enough: the syntax obtained that way was still too strong. When proving properties of Cartegena logic, we noticed that on many occasions we had two different partitions, and we had to take a common partition. The notion of downwards correctness encapsulates when truth of a formula is preserved when a partition is replaced with a finer one. Explicitly, we needed Lemmas \ref{correct} and \ref{lem}, which, in fact, are the motivation for the concepts of upwards/downwards correctness. These lemmas make transparent when replacing a partition with a finer one preserves the truth of a formula.

We will now define upwards and downwards correct formulas.

\begin{de}[Good function]
    A function \[
p:\Pow(\theta)\to \Pow(\Pow(\theta))
\] is called \textit{good} if
\begin{enumerate}
    \item $p(A)\subseteq\Pow(A)$,
    \item $A'\subseteq A \implies  p(A')\subseteq p(A)$,
\end{enumerate}
for all $A',A\subseteq\theta$.
\end{de} 

\begin{ex} The following are the most important good functions:
\begin{align*}
    &A\mapsto\emptyset.\\
    & A\mapsto[A]^n.\\
    & A\mapsto[A]^\theta.\\
    & A\mapsto\Pow(A).
\end{align*}
Here $n$ is a natural number and $\theta$ is a cardinal.
\end{ex}

A formula being ``upwards correct with respect to $X$'' can be understood as an analogue of a ``positive occurrence of a subformula''. Instead of a subformula, we have the variable $X$. On the other hand, ``downwards correctness'' is then an analogue for ``negative occurrence''.

In the following inductive definition, the crucial steps are the two last ones: they describe how upwards correctness behaves with respect to infinitary disjunctions and downwards correctness behaves with respect to infinitary conjunctions (along good functions).

\begin{de}[Upwards/downwards correctness] Let $\kappa$ be a cardinal, let $X\in\BVar_\kappa^\theta$ be a Boolean variable and let $\phi$ be an $\LL^{\Bool}_{\kappa\kappa}$-formula.
\begin{enumerate}
    \item If neither $X$ nor any of $x_i^X$, $i\in\theta$, occur free in $\phi$, then $\phi$ is upwards and downwards correct for $X$.
    \item If $\phi$ is upwards correct for $X$, then $\neg\phi$ is downwards correct for $X$.
    \item If $\phi$ is downwards correct for $X$, then $\neg\phi$ is upwards correct for $X$.
    \item If $\phi$ is upwards correct for $X$, then $\exists x\phi$  and $\forall x\phi$ are upwards correct for $X$.
    \item If $\phi$ is downwards correct for $X$, then $\exists x\phi$  and $\forall x\phi$ are downwards correct for $X$.
    \item If $\Phi$ is a set of formulas such that every $\phi\in\Phi$ is upwards correct for $X$, then $\bigvee\Phi$ and $\bigwedge\Phi$ are upwards correct for $X$.
    \item If $\Phi$ is a set of formulas such that every $\phi\in\Phi$ is downwards correct for $X$, then $\bigvee\Phi$ and $\bigwedge\Phi$ are downwards correct for $X$.
    \item\label{88} If
    \begin{enumerate}
        \item $\{\phi_u:u\subseteq\theta\}$ are $\LL^{\Bool}_{\kappa\kappa}$-formulas,
        \item $X$ does not occur in any $\phi_u$,
        \item\label{cc} $\FreeVar(\phi_u)\cap\{x_i^X:i\in\theta\}\subseteq\{x_i^X:i\in u\}$ for every $u\subseteq\theta$,
        \item $p$ is a good function on $\Pow(\theta)$,
    \end{enumerate}
    then
    \[
        \bigdvee_{u\in p(X)}\phi_u,
    \]
    is upwards correct with respect to $X$.
    \item\label{99} If
    \begin{enumerate}
        \item $\{\phi_u:u\subseteq\theta\}$ are $\LL^{\Bool}_{\kappa\kappa}$-formulas,
        \item $X$ does not occur in any $\phi_u$,
        \item\label{ccc} $\FreeVar(\phi_u)\cap\x^X_\theta\subseteq\x_u^X$ for every $u\subseteq\theta$,
        \item $p$ is a good function on $\Pow(\theta)$,
    \end{enumerate}
    then
    \[
        \bigdwedge_{u\in p(X)}\phi_u,
    \]
    is downwards correct with respect to $X$.
\end{enumerate}
    
\end{de}

\begin{rmk}
    An atomic formula is upwards or downwards correct for a Boolean variable $X\in\BVar^\theta_\kappa$ if and only if none of the variables $\{x_i^X:i\in\theta\}$ occurs in it. This is to prevent formulas such as
    \[
    \bigwedge_{n\in\omega}x_{n+1}<x_n
    \]
    from being allowed to be closed into a sentence. Clauses \ref{cc} and \ref{ccc} are crucial here.
\end{rmk}

\begin{notation}
    Let $\phi=\phi(\x,\X)$ and $\psi=\psi(\x,\X)$ be $\LL^{\Bool}_{\kappa\kappa}$-formulas. We write
    \[
    \phi\models\psi
    \]
    if for all structures $\MM$, all tuples $\a\in\MM^{\x}$ and all valuations $\A$ for $\X$, 
    \[
    \MM\models\phi(\a,\A)\implies\MM\models\psi(\a,\A).
    \]
\end{notation}

The next two lemmas explain whether the truth of a formula is preserved when a value of a Boolean variable is replaced by another (Lemma \ref{correct}), or when a partition is replaced by a finer or a coarser one (Lemma \ref{lem}). These lemmas play a key role in the proof of game-syntax correspondence (Theorem \ref{gamesyntax}). They also explain the names ``upwards/downwards correctness''.

\begin{lem}\label{correct}
    Let $X$ be a Boolean variable ranging over subsets of some cardinal $\theta$ and let $\phi$ be a formula.
    \begin{enumerate}
        
    \item If $\phi$ is upwards correct for $X$, then for all $A'\subseteq A\subseteq\theta$, 
    \[
    \phi[A'/X]\models\phi[A/X].
    \]
    \item If $\phi$ is downwards correct for $X$, then for all $A'\subseteq A\subseteq\theta$, 
    \[
    \phi[A/X]\models\phi[A'/X].
    \]
    \end{enumerate}
\end{lem}
\begin{proof} 
    Follows inductively by applying monotonicity of good functions.
\end{proof}

Lemma \ref{correct} immediately gives the following:

\begin{lem}\label{lem} 
Let $X\in\BVar_\kappa^\theta$ and let $\WW$ and $\VV$ be partitions of $\theta$ such that $\WW$ refines $\VV$. Let $\phi$ be an $\LL^{\Bool}_{\kappa\kappa}$-formula.
\begin{enumerate}
    \item If $\phi$ is downwards correct for $X$, then 
        \[
        \bigdwedge_{X\in\VV}\phi\models\bigdwedge_{X\in\WW}\phi.
        \]
    \item If $\phi$ is upwards correct for $X$, then 
    \[
        \bigdvee_{X\in\WW}\phi\models\bigdvee_{X\in\VV}\phi.
        \]
\end{enumerate}
\end{lem}

Finally we are ready to define Cartagena formulas.

\subsubsection{Cartagena formulas}\spa 

\vv

First recall that for $f:\theta\to\omega$,
\[
\WW_f:=\{f^{-1}\{n\}:n\in\omega\}.
\]

\begin{de}[Cartagena logic $\LL^c_\kappa$]\label{syntax}  Let $\lambda$ be an infinite regular cardinal. The $\LL^c_\lambda$-formulas are defined as follows.
    \begin{enumerate}
        \item\label{1'} Every atomic formula with variables in $\Var_\lambda$ is an $\LL^c_\lambda$-formula.
        \item\label{2'} If $\phi$ is an $\LL^c_\lambda$-formula, then $\neg\phi$ is an $\LL^c_\lambda$-formula.
        \item\label{3'} If $\Phi$ is a set of size $<\lambda$ of $\LL^c_\lambda$-formulas with \[
        |\bigcup_{\phi\in\Phi}\FreeVar(\phi)\cup\FreeBVar(\phi)|<\lambda,
        \] then $\bigwedge\Phi$ and $\bigvee\Phi$ are $\LL^c_\lambda$-formulas.
        \item If $\{\phi_u:u\subseteq\theta\}$ are $\LL^c_\lambda$-formulas, $X\in\BVar_\lambda^\theta$ and
        \[
        \phi=\bigdvee_{u\in p(X)}\phi_u,
        \]
        is upwards correct with respect to $X$, then $\phi$ is an $\LL^c_\lambda$-formula.
        \item If $\{\phi_u:u\subseteq\theta\}$ are $\LL^c_\lambda$-formulas, $X\in\BVar_\lambda^\theta$ and
        \[
        \phi=\bigdwedge_{u\in p(X)}\phi_u
        \]
        is downwards correct with respect to $X$, then $\phi$ is an $\LL^c_\lambda$-formula.
        \item\label{6'} If $\phi$ is an $\LL^c_\lambda$-formula, then so are $\exists x\phi$ and $\forall x\phi$. 
        
        \item\label{7'} If $\phi$ is an $\LL^c_\lambda$-formula such that \begin{enumerate}
            \item $\phi$ is downwards correct for $X$,
            \item $\FreeVar(\phi[A/X])\cap\{x^X_i:i\in\theta\}\subseteq\{x_i^X:i\in A\}$ for all $A\subseteq\theta$,
        \end{enumerate} then the following are $\LL^c_\lambda$-formulas:
        \begin{align*}
            & \exists\x^X_\theta\bigvee_{f:\theta\to\omega}\bigdwedge_{X\in\WW_f}\phi,\\
            & \forall\x^X_\theta\bigvee_{f:\theta\to\omega}\bigdwedge_{X\in\WW_f}\phi.
        \end{align*}
        \item\label{8'} If $\phi$ is an $\LL^c_\lambda$-formula such that \begin{enumerate}
            \item $\phi$ is upwards correct for $X$,
            \item $\FreeVar(\phi[A/X])\cap\{x^X_i:i\in\theta\}\subseteq\{x_i^X:i\in A\}$ for all $A\subseteq\theta$,
        \end{enumerate} then the following are $\LL^c_\lambda$-formulas:
        \begin{align*}
            & \exists\x^X_\theta\bigwedge_{f:\theta\to\omega}\bigdvee_{X\in\WW_f}\phi,\\
            & \forall\x^X_\theta\bigwedge_{f:\theta\to\omega}\bigdvee_{X\in\WW_f}\phi.
        \end{align*}
        
    \end{enumerate}
    For singular $\kappa$, we define
    \[
    \LL^c_\kappa:=\bigcup_{\lambda<\kappa}\LL^c_{\lambda^+}.
    \]
    Whenever $\lambda$ is clear from context, we refer to $\LL^c_\lambda$-formulas as \textit{Cartagena formulas}.
\end{de}

\begin{rmk}
    As the logic is closed under negation, some of the clauses are redundant. 
\end{rmk}

\begin{rmk}
    The notion of subformula is different in $\LL^c_\lambda$ from the one in $\LL_{\lambda\lambda}$. Indeed, we do not ``break'' Cartagena quantifiers: the set of subformulas of the formula
    \[
    \exists\x^X_\theta\bigvee_{f:\theta\to\omega}\bigdwedge_{X\in\WW_f}\phi
    \]
    is the set $\{\phi\}$ together with the set of subformulas of $\phi$.
\end{rmk}

\begin{notation}
    In the absence of a risk of confusion, we omit the upper index $X$ from the variables $\x_\theta^X$ and write
    \[
    \exists\x_\theta\bigvee_{f:\theta\to\omega}\bigdwedge_{X\in\WW_f}\phi,
    \]
    thereby assuming that for each $A\subseteq\theta$,
    \[
    \FreeVar(\phi[A/X])\cap\{x_i:i\in\theta\}\subseteq\{x_i:i\in A\}.
    \]
    Similarly for a pair $\y_\theta$ and $Y$, etc. For instance, in writing 
    \[
    \exists\x_\theta\bigvee_{f:\theta\to\omega}\bigdwedge_{X\in\WW_f}\left(\forall\y_\theta\bigwedge_{g:\theta\to\omega}\bigdvee_{Y\in\WW_g}\phi(\x_\theta,\y_\theta,X,Y)\right),
    \]
    it is implicit that $\x_\theta^X$ and $\y_\theta^Y$ are identified with $\x_\theta$ and $\y_\theta$, respectively, and that ${\{x_i:i\in\theta\}\cap\{y_i:i\in\theta\}=\emptyset}$. Similarly for other Cartagena quantifiers.
\end{notation}

\begin{rmk}
    The semantics for Cartagena formulas follow from \ref{semantics}. In particular,
    \[
    \MM\models\exists\x_\theta\bigvee_{f:\theta\to\omega}\bigdwedge_{X\in\WW_f}\phi(\x_\theta,X)
    \]
    holds if and only if there are $\a_\theta\in\MM^\theta$ and $f:\theta\to\omega$ such that for all $A\in\WW_f$, 
    \[
    \MM\models\phi(\a_A,A).
    \]
    The other clauses are defined similarly.
\end{rmk}

\begin{de}
    The \textit{quantifier rank} $\qr(\phi)$ of a Cartagena formula $\phi$ is defined as follows:

\begin{align*}
    &\qr(\phi):=0,\quad\text{for atomic }\phi,\\
    &\qr(\neg\phi):=\qr(\phi),\\
    &\qr\left(\bigwedge\Phi\right):=\sup_{\phi\in\Phi}(\qr(\phi)),\\
    &\qr(\exists x\phi):=\qr(\phi)+1,\\
    &\qr\left(\bigdwedge_{u\in p(X)}\phi_u\right):=\sup_{u\subseteq\theta}(\qr(\phi_u)),\\
    &\qr\left(\exists\x_\theta\bigvee_{f:\theta\to\omega}\bigdwedge_{X\in\WW_f}\phi\right):=\qr(\phi)+1,\\
    &\qr\left(\forall\x_\theta\bigvee_{f:\theta\to\omega}\bigdwedge_{X\in\WW_f}\phi\right):=\qr(\phi)+1.
\end{align*}
\end{de}

Notice the quantifier rank is defined with respect to standard variables $\Var$. A Boolean variable does not affect the quantifier rank. Furthermore, in general
\[
\qr(\phi(\x,\X))\geq\qr(\phi(\x,\A)),
\]
where $\A$ is a valuation of $\X$. A strict inequality is possible.

We are finally ready to prove that Cartagena syntax corresponds to Cartagena game.

\subsection{Game-syntax equivalence}\spa 

\vv

\begin{de}
    We write
    \[
    \MM\equiv^\beta_\kappa\NN
    \]
    if $\MM$ and $\NN$ agree on $\LL^c_\kappa$-sentences of quantifier rank $\leq\beta$.
\end{de}

\begin{rmk}\label{counting} If $\kappa$ is an uncountable cardinal such that $\kappa=\beth_\kappa$, then there are $<\kappa$ many $\LL^{\Bool}_{\kappa\kappa}$-formulas up to equivalence of quantifier rank $<\beta$, in $<\kappa$ many free variables.

In detail: For each regular $\lambda<\kappa$, there are at most $\beth_{2\cdot\beta}(\eta+\lambda)$ many $\LL_{\lambda\lambda}$-formulas up to equivalence of quantifier rank $<\beta$ in $\leq\eta$ many free variables, for all ordinals $\beta<\lambda$, following Benda \cite{benda1969reduced}. The projection map $\pi_{\A}$  from Proposition \ref{collapsing} is ${<\lambda}$-to-one, because each formula is a well founded tree with levels of size $<\lambda$. The number of $\LL^{\Bool}_{\lambda\lambda}$-formulas of quantifier rank $<\beta$ in $\leq\eta$ many free standard variables and $\leq\lambda$ many free Boolean variables $\X\subseteq\BVar_\lambda$ is therefore at most
    \[
    \left|\bigcup_{\A\text{ valuation for }\X}\pi_{\A}^{-1}\left[\LL^{<\beta}_{\lambda\lambda}\right]\right|\leq 2^\lambda\cdot\lambda\cdot\beth_{2\cdot\beta}(\eta+\lambda)\leq \beth_{2\cdot\beta}(\eta+\lambda),
    \]
which is strictly less than $\kappa$ if $\kappa=\beth_\kappa$.
\end{rmk}

Recall that for functions $f,g:\theta\to\omega$, we write
\[
    f\leq g\quad:\iff\quad \WW_f\text{ refines }\WW_g.
\]

\begin{thm}[Game-syntax correspondence]\label{gamesyntax} Assume that $\kappa$ is an infinite cardinal such that $\kappa=\beth_\kappa$.
    The following are equivalent:
    \begin{enumerate}
        \item\label{gseka} $\MM\equiv^\beta_\kappa\NN$.
        \item\label{gstoka} Player $\2$ has a winning strategy in the game $\G^\beta_\kappa(\MM,\NN)$.
    \end{enumerate}
\end{thm}
\begin{proof}\spa 

    \noindent \textit{(\ref{gseka}) $\To$ (\ref{gstoka})}: The proof is by induction on $\beta$. We assume that the claim holds for all $\beta'<\beta$ and suppose that $\MM\equiv^\beta_\kappa\NN$. We claim that player $\2$ has a winning strategy in the game $\G^\beta_\kappa(\MM,\NN)$. 
        
        It suffices to show that player $\2$ has a strategy such that each state $(\alpha,\pi)$ reached in the game satisfies:
        \[
        \textit{The partial isomorphism }\pi\textit{ is }\alpha\textit{-elementary}. \tag{$*$}\label{star}
        \]
        Here $\alpha$\emph{-elementary} means that $\pi$ preserves formulas of quantifier rank $\leq\alpha$ with any values for Boolean varibles.

        The starting state $(\beta,\emptyset)$ clearly satisfies the condition \eqref{star}.
        
        Suppose that the game is at state $(\alpha,\pi)$, and $(\alpha,\pi)$ satisfies \eqref{star}. We show that whichever is the next move of player $\1$, player $\2$ can play in such a way that the next state satisfies ($*$).  
        
        To this end, at state $(\alpha,\pi)$, suppose that player $\1$ plays an ordinal $\alpha'<\alpha$ and a tuple $\a_\theta$. By symmetry, we may assume that $\a_\theta\in\MM^\theta$. First, we enumerate the domain of $\pi$ as $\bar{d}=\dom(\pi)$ and denote $\pi(\bar{d}):=(\pi(d_i))_i$.
        For every $u\subseteq\theta$, let
            \[
            \phi_u(\x_u,\y):=\bigwedge\{\psi(\x_u,\y):\qr(\psi)\leq\alpha'\text{ and }\MM\models\psi(\a_u,\bar{d})\}.
            \]
            This is a conjunction of size $<\kappa$, by Remark \ref{counting}.  Let
            \[
            \phi(\x_\theta,\y,X):=\bigdwedge_{u\in \Pow(X)}\phi_u(\x_u,\y).
            \]
            This is a Cartagena formula - downwards correct for the Boolean variable $X$. In this formula, the index function $p$ is the full powerset function, $X\mapsto p(X):=\Pow(X)$, which is a good function.
            Furthermore, we have
            \[
            \MM\models\exists\x_\theta\bigvee_{f:\theta\to\omega}\bigdwedge_{X\in\WW_f}\phi(\x_\theta,\bar{d},X).
            \]
            This formula $\exists\x_\theta\bigvee_{f:\theta\to\omega}\bigdwedge_{X\in\WW_f}\phi(\x_\theta,\y,X)$ is a Cartagena formula, since the formula $\phi(\x_\theta,\y,X)$ is downwards correct for $X$. It has quantifier rank $\alpha'+1\leq\alpha$. By the assumption that $\pi$ is $\alpha$-elementary, we have 
            \[
            \NN\models\exists\x_\theta\bigvee_{f:\theta\to\omega}\bigdwedge_{X\in\WW_f}\phi(\x_\theta,\pi(\bar{d}),X).
            \]
            Let $\b_\theta\in\NN^\theta$ and $f:\theta\to\omega$ be such that 
            \[
            \NN\models \bigdwedge_{X\in\WW_f}\phi(\b_\theta, \pi(\bar{d}),X).
            \] 
            Now we let player $\2$ play $\b_\theta$ and $f$. For every $A\in\WW_f$, the map
            \[
            \pi_A:=\pi\cup\{(a_i,b_i):i\in A\}
            \]
            is $\alpha'$-elementary. Thus, whichever piece $A\in\WW_f$ player $\1$ chooses to play, the next state $(\alpha',\pi_A)$ satisfies \eqref{star}.
            
\vv

\noindent \textit{(\ref{gstoka}) $\To$ (\ref{gseka})}: The proof is by induction on $\beta$. Our induction hypothesis is slightly stronger. Write
\begin{align*}
    \mathsf{IH}_\alpha\spa :=\spa & \forall\a \in\MM^{<\kappa}\spa \forall\b\in\NN^{<\kappa}\quad \2\uparrow\G_\kappa^\alpha((\MM,\a),(\NN,\b)) \to (\MM,\a)\equiv^\alpha_\kappa(\NN,\b).
\end{align*}
We assume that $\IH_\alpha$ holds for every $\alpha<\beta$ and show that $\IH_\beta$ holds.

For simplicity, we ignore the parameters $\a$ and $\b$. Suppose that player $\2$ has a winning strategy in $\G_\kappa^\beta(\MM,\NN)$. We argue by induction on the complexity of a formula $\psi$ of quantifier rank $\leq\beta$ that $\MM$ and $\NN$ agree on $\psi$. The claim is clear when $\psi$ is atomic. For the inductive steps, we check only the Cartagena clauses. 

\vv

\noindent \textbf{Case 1:} $\psi\mbox{ is } \forall\x_\theta\bigvee_{f:\theta\to\omega}\bigdwedge_{X\in\WW_f}\phi(\x_\theta,X)$.

\vv

Assuming $\MM\models\psi$, we show that $\NN\models\psi$.
Let $\b_\theta\in\NN^\theta$ be arbitrary. It suffices to find $f:\theta\to\omega$ such that for all $A\in\WW_f$, 
\[
\NN\models \phi(\b_A,A).
\] 
Let player $\1$ play the tuple $\b_\theta$ and the ordinal $\alpha:=\qr(\phi)$.
The winning strategy of player $\2$ gives a tuple $\a_\theta\in\MM^\theta$ and a function ${h:\theta\to\omega}$.
We apply the fact that $\MM\models\psi$ to the tuple $\a_\theta$ and obtain $g:\theta\to\omega$ such that for all $A\in\WW_g$,
\[
\MM\models\phi(\a_A,A).
\]
Let $f\leq g,h$.
As $\phi$ is downwards correct for $(X,\x_\theta)$ and $f\leq g$, we have
\[
\MM\models\phi(\a_A, A)
\]
for every $A\in\WW_f$.
This follows from Lemma \ref{correct}.
By Lemma \ref{strategy}, player $\2$ can play the function $f$ instead of $h$, and preserve her winning strategy - she thus has a winning strategy in the further game $\G^\alpha_\kappa((\MM,\a_A),(\NN,\b_A))$, for every $A\in\WW_f$. By $\IH_\alpha$, for each $A\in\WW_f$, 
\[
(\MM,\a_A)\equiv^\alpha_\kappa(\NN,\b_A).
\]
Hence for all $A\in\WW_f$,
\[
\NN\models \phi(\b_A,A),
\]
as wanted. Thus $\NN\models\psi$.

\vv

\noindent\textbf{Case 2:} $\psi\mbox{ is }\exists\x_\theta\bigvee_{f:\theta\to\omega}\bigdwedge_{X\in\WW_f}\phi(\x_\theta,X)$.

\vv

This proof is very similar to the one above, with an exception that we make player $\1$ play first in $\MM$ a witness for the existential quantifier. The response of player $\2$, a tuple in $\NN$, will be a witness for the existential quantifier in $\NN$. We can again take a common refinement of the function given by the disjunction true in $\MM$ and the one given by player $\2$ - and this function will work.

Suppose that $\MM\models \psi$. We show that $\NN\models\psi$. It suffices to find $\b_\theta\in\NN^\theta$ and $f:\theta\to\omega$ such that
    \[
    \NN\models\phi(\b_A,A)
    \]
    for all $A\in\WW_f$.
Let $\a_\theta\in\MM^\theta$ and $g:\theta\to\omega$ be such that for all $A\in\WW_g$, 
    \[
    \MM\models\phi(\a_A,A).
    \]
Let player $\1$ play the tuple $\a_\theta$ and the ordinal $\alpha:=\qr(\phi)$.
The winning strategy of player $2$ gives a tuple $\b_\theta\in\NN^\theta$ and a function ${h:\theta\to\omega}$.
Let $f\leq g,h$.
The rest of the proof is verbatim the same as the proof of Case 1. In the end we obtain that
    \[
    \NN\models\phi(\b_A,A)
    \]
    for all $A\in\WW_f$, as wanted, which shows that indeed $\NN\models\psi$.

\vv 

This ends the proof.

\end{proof}

\section{Model theory of Cartagena logic}\label{mtocl}

In this section we give three important model-theoretic properties of Cartagena logic: closure under unions of countable elementary chains (Theorem \ref{union}), Löwenheim-Skolem-Tarski Theorem (Theorem \ref{LST}) and a strong form of undefinability of well order (Theorem \ref{sudwo}). The first two of these are strengthenings of known analogous properties of $\LL^1_\kappa$. The generative syntax of Cartagena logic was the key that made possible to formulate and prove the stronger versions. The third property, strong undefinability of well order, holds for $\LL^1_\kappa$ as such.

We start by giving the definition of a fragment, which will be used in all three proofs.

Throughout the section, let $\kappa$ be an uncountable cardinal with $\kappa=\beth_\kappa$.

\begin{de}
    A \textit{fragment (of Cartagena logic)} is a set $\LL$ of $\LL^c_\kappa$-formulas in a fixed signature $\tau$ containing all first order $\tau$-formulas, closed under first order operations, taking subformulas, substituting a term for a standard variable and substituting any valuation for a Boolean variable.
\end{de}

\begin{de}\label{fragelem}
    For a fragment $\LL$, we write
    \[
    \NN\elem_\LL\MM
    \]
    and say that $\NN$ is an \textit{$\LL$-elementary submodel} of $\MM$, if $\NN\subseteq\MM$, and for every $\phi(\x,\X)\in\LL$: for all $\a\in\NN^{\x}$ and for every valuation $\A$ of $\X$,
    \[
    \NN\models\phi(\a,\A)\quad\iff\quad\MM\models\phi(\a,\A).
    \]    A chain $(\MM_i)_{i\in\delta}$ is \textit{$\LL$-elementary} if $\MM_j\elem_\LL\MM_i$ for all $j<i$.
\end{de}

\subsection{Union Lemma}\spa 

\vv

The Union Lemma \ref{union} of Cartagena logic, that states that Cartagena logic is closed under unions of countable elementary chains, is not known to hold for $\LL^1_\kappa$ as stated here\footnote{It is, however, possible to borrow a stronger elementary substructure relation from $\LL_{\kappa\kappa}$ and prove a weaker version of the Union Lemma for $\LL^1_\kappa$. See the explicit statement in \cite{shelah2012nice}.}. The problem with the logic $\LL^1_\kappa$, lacking a simple notion of a formula with free variables, is that it does not have the same definition of elementary substructure as logics usually do, rendering the notion of elementary chain more difficult to treat. 

\begin{thm}[Union Lemma]\label{union} Let $\LL$ be a fragment. If $(\MM_n)_{n\in\omega}$ is an $\elem_\LL$-chain, then for every $n\in\omega$,
\[
\MM_n\elem_\LL\bigcup_{n\in\omega}\MM_n.
\]
    
\end{thm}
\begin{proof}
    Let 
    \[
    \MM:=\bigcup_n\MM_n.
    \]
    It is enough to show that $\MM_0\elem_\LL\MM$. We show by induction on the complexity of $\psi(\x,\X)$ that for each valuation $\A$ for $\X$,
    \[
    \forall\a\in\MM_0^{\x}:\quad \MM_0\models\psi(\a,\A)\quad\iff\quad\MM\models\psi(\a,\A).
    \]
    We only deal with the clauses that start with a Cartagena quantifier. The other clauses are standard.

    \vv 

    \noindent\textbf{Case 1:} $\psi\mbox{ is }\exists\x_\theta\bigvee_{f:\theta\to\omega}\bigdwedge_{X\in\WW_f}\phi(\x_\theta,\y,X,\Y)$.

    \vv

    For simplicity, we suppose that $\y=\emptyset=\Y$ and that the only free variables in $\phi$ are $\x_\theta$ and $X$. Note that the formula $\phi(\x_\theta,X)$ is downwards correct for the variable $X$, by assumption. 
    
    \vv 
    
    We first show
    \[
    \MM\models\psi\implies\MM_0\models\psi.
    \]
    Let $\a_\theta\in\MM^\theta$ and $f:\theta\to\omega$ be such that for every $A\in\WW_f$,
    \[
    \MM\models\phi(\a_A,A).
    \]
    Up to refining $f$, we may assume that for every $A\in\WW_f$ there is $n_A$ such that $\a_A\subseteq\MM_{n_A}$. Refining $f$ is allowed as $\phi$ is downwards correct for $X$, by Lemma \ref{correct}.
    By induction hypothesis, for every $A\in\WW_f$, 
    \[
    \MM_{n_A}\models \phi(\a_A,A), 
    \]
    which implies
    \[
    \MM_{n_A}\models\exists\y_A\bigvee_{g:A\to\omega}\bigdwedge_{Y\in \WW_g}\phi(\y_A, Y).
    \]
    This is a correct Cartagena formula up to re-enumerating $A$ with its cardinality.
    As $\MM_0\elem_\LL\MM_{n_A}$ for each $A\in\WW_f$, we have
    \[
    \MM_0\models\exists\y_A\bigvee_{g:A\to\omega}\bigdwedge_{Y\in \WW_g}\phi(\y_A,Y).
    \]
    For each $A\in\WW_f$ pick $\b_A\in\MM_0^A$ and $g_A:A\to\omega$ such that
    \[
    \MM_0\models\bigdwedge_{Y\in\WW_g}\phi(\b_Y,Y).
    \]
    As $\WW_f$ is a partition of $\theta$, we can glue the tuples $(\b_A)_{A\in\WW_f}$ and the functions $(g_A:A\to\omega)_{A\in\WW_f}$ together: we have $\b_\theta:=(b_i)_{i\in\bigcup\WW_f}$, and
    \begin{align*}
        g:\theta\to\omega, \quad g(i):=g_A(i),\text{ for the unique }A\in\WW_f\text{ with }i\in A.
    \end{align*}
    These two witness
    \[
    \MM_0\models\exists\y_\theta\bigvee_{g:\theta\to\omega}\bigdwedge_{Y\in\WW_g}\phi(\y_\theta,Y),
    \]
    which is as wanted.

    The direction
    \[
    \MM_0\models\psi\implies\MM\models\psi
    \]
    is an immediate application of the induction hypothesis.

    \vv 

    \noindent\textbf{Case 2:} $\psi\mbox{ is }\forall\x_\theta\bigvee_{f:\theta\to\omega}\bigdwedge_{X\in\WW_f}\phi(\x_\theta,\y,X,\Y)$. 

    \vv 

    Again, for simplicity, we suppose that $\y=\emptyset=\Y$ and that the only free variables in $\phi$ are $\x_\theta$ and $X$.
    Again, the formula $\phi$ is downwards correct for $X$.
    
\vv 

The direction
        \[
        \MM\models\psi\implies\MM_0\models\psi.
        \]
        is an immediate application of the induction hypothesis.

We show 
        \[
        \MM_0\models\psi\implies\MM\models\psi.
        \]
        Let $\a_\theta\in\MM^\theta$. We find $f:\theta\to\omega$ such that for every $A\in\WW_f$,
        \[
        \MM\models\phi(\a_A,A).
        \]
        Let 
        \[
        g:\theta\to\omega,\quad g(i):=\text{the least }n\text{ such that }a_i\in\MM_n.
        \]
        For each $n$, let $B_n:=g^{-1}\{n\}$. Notice that
        \[
        \a_{B_n}\in\MM_n^{B_n}.
        \]
        As $\MM_n\models\psi$ for each $n$, we also have
        \[
        \MM_n\models\forall\x_{B_n}\bigvee_{f:B_n\to\omega}\bigdwedge_{Y\in\WW_f}\phi(\x_{B_n},Y),
        \]
        for each $n$. These are Cartagena formulas, up to enumerating $B_n$ by its cardinality. We apply the universal quantifiers $\forall\x_{B_n}$ to the tuples $\a_{B_n}$. For each $n$, let $g_n:B_n\to\omega$ be such that for all $A\in\WW_{g_n}$
        \[
        \MM_n\models\phi(\a_A,A).
        \]
        By induction hypothesis for every $n\in\omega$ and every $A\in\WW_{g_n}$,
        \[
        \MM\models\phi(\a_A,A).
        \]
        Note that $g_n\leq g\rest B_n$ for each $n$. As $(B_n)_{n\in\omega}$ is a partition of $\theta$, we can glue the functions $(g_n:B_n\to\omega)_{n\in\omega}$ together:
        \[
        f:\theta\to\omega,\quad f(i):=g_n(i),\text{ for the unique }n\text{ such that }i\in B_n.
        \]
        Now $f:\theta\to\omega$ is such that for every $A\in\WW_f$,
        \[
        \MM\models\phi(\a_A,A),
        \]
        which is as wanted.

\vv 

    This ends the proof.
    
\end{proof}

\subsection{Löwenheim-Skolem-Tarski Theorem}\spa 

\vv 

In this section, we prove a Löwenheim-Skolem-Tarski Theorem for Cartagena logic. The proof is similar to that of first order logic - using Skolem functions.

\begin{de} Let $\phi$ be an $\LL^c_\kappa$-formula. The set
    \[
\LL_\phi:=\bigcap\{\LL:\LL\textit{ is  a fragment and }\phi\in\LL\}
\]
is called the \textit{fragment below $\phi$}.
\end{de}

The fragment below $\phi$ is a fragment of size $<\kappa$ such that $\phi\in\LL_\phi$, for every $\LL^c_\kappa$-formula $\phi$. This follows from Remark \ref{counting}.

\begin{thm}[Löwenheim-Skolem-Tarski]\label{LST} Let $\LL$ be a fragment of Cartagena logic.
For every model $\MM$ and every $E\subseteq\MM$ there is a model $\NN$ such that
    \[
    E\subseteq\NN\elem_\LL\MM\quad\text{and}\quad |\NN|\leq \min\{|E|+\kappa,(|E|+|\LL|)^+\}.
    \]
\end{thm}
\begin{proof} Let $\LL$ be a fragment, let $\MM$ be a model and choose a subset $E\subseteq\MM$. We will find functions $f_\phi$ and $g_\psi$ that will act as Skolem functions and close $E$ under these.

Let $\BB$ be the set of Boolean variables occurring in the fragment $\LL$. For a tuple of Boolean variables $\Y\subseteq\BB$, write
\[
\Val(\Y):=\{\B:\B\text{ is a valuation for }\Y\}.
\]
Notice that $|\Val(\BB)|\leq|\LL|$.
For every formula $\phi\in\LL$ which is upwards correct for $X$ and for every formula $\psi$ which is downwards correct for $X$ such that
\begin{align*}
    & \exists\x_\theta\bigvee_{f:\theta\to\omega}\bigdwedge_{X\in\WW_f}\phi(\x_\theta,\y,X,\Y)\in\LL,\\
    & \exists\x_\theta\bigwedge_{f:\theta\to\omega}\bigdvee_{X\in\WW_f}\psi(\x_\theta,\y,X,\Y)\in\LL,
\end{align*}
there are Skolem functions
\begin{align*}
    & f_\phi:\MM^{\y}\times\mathsf{Val}(\Y)\to\MM^\theta,\\
    & g_\psi:\MM^{\y}\times\mathsf{Val}(\Y)\to\MM^\theta,
\end{align*}
such that for all $(\b,\B)\in \MM^{\y}\times\mathsf{Val}(\Y)$, 
\begin{align*}
    & \MM\models \exists\x_\theta\bigvee_{f:\theta\to\omega}\bigdwedge_{X\in\WW_f}\phi(\x_\theta,\b,X,\B)\to\bigvee_{f:\theta\to\omega}\bigdwedge_{X\in\WW_f}\phi(f_\phi(\b, \B),\b,X,\B),\\
    & \MM\models \exists\x_\theta\bigwedge_{f:\theta\to\omega}\bigdvee_{X\in\WW_f}\psi(\x_\theta,\b,X,\B)\to\bigwedge_{f:\theta\to\omega}\bigdvee_{X\in\WW_f}\psi(g_\psi(\b, \B),\b,X,\B).
\end{align*}
    
    We let $\NN$ to be the closure of $E\times\Val(\BB)$ under the functions $f_\phi$, $g_\psi$. In other words, $\NN$ is obtained by first denoting
    \begin{align*}
        \FF:=& \{f_\phi:\exists\x_\theta\bigvee_f\bigdwedge_X\phi(\x,\y,X,\Y)\in\LL\}\cup\{g_\psi:\exists\x_\theta\bigwedge_f\bigdvee_X\psi(\x,\y,X,\Y)\in\LL\}.
    \end{align*}
    and then closing recursively:
    \begin{itemize}
        \item $\NN_0:=E$.
        \item $\NN_{\alpha+1}:=\NN_n\cup\bigcup_{f\in\FF}f[\NN_\alpha\times\Val(\BB)]$,
        \item $\NN_\alpha:=\bigcup_{\beta<\alpha}\NN_\beta$, for limit ordinal $\alpha$.
    \end{itemize}
    Finally, $\NN:=\NN_{\min\{|E|+\kappa,(|E|+|\LL|)^+\}}$. It now follows that $\NN$ is as wanted: $E\subseteq\NN$, $\NN$ is a substructure of $\MM$ of size $\min\{|E|+\kappa,(|E|+|\LL|)^+\}$ and $\NN\elem_\LL\MM$.
\end{proof}

Given a Cartagena sentence $\phi$, the fragment $\LL_\phi$ has size $<\kappa$, and therefore by Theorem $\ref{LST}$ for every model $\MM$ there is a model $\NN\elem_{\LL_\phi}\MM$ of size $<\kappa$. In particular, such $\NN$ agrees on $\phi$ with $\MM$. We obtain:

\begin{cor}
    Every Cartagena sentence that has a model has a model of size $<\kappa$.
\end{cor}

\subsection{Strong undefinability of well order}\label{sectsudwo}\spa 

\vv

This section provides a proof of the fact that Cartagena logic $\LL^c_\kappa$ has \textit{strong well ordering number $\kappa$}. This property is a slight modification of a property due to Shelah \cite{shelah2012nice}. It combines a weak form of compactness and a Löwenheim-Skolem type of property, and it is important because the logic $\LL^1_\kappa$ is maximal with respect to it (Theorem 3.4 in \cite{shelah2012nice}).

The setup is very general, and concerns not only Cartagena logic $\LL^c_\kappa$, but also $\LL^1_\kappa$, and any abstract logic satisfying the assumptions of the theorem. We include the proof for the sake of completeness. A simplification of the proof in the case of $\LL_{\kappa\omega}$ can be found in \cite{väänänen_2011}.

We start with the relevant definitions.

\begin{de}\label{sudwodef} Let $P$ be a unary predicate and let $<$ be a binary predicate symbol.
    \begin{enumerate}
        \item The \textit{well ordering number} of a logic $\LL$, if it exists, is the least ordinal $\alpha$ such that for any $\LL$-sentence $\phi$ in a signature containing $P$ and $<$ the following holds: if $\phi$ has a model $\MM$ such that
        \[
        (P^\MM,<^\MM)\cong(\alpha,\in),
        \]
        then $\phi$ also has a model $\NN$ such that $(P^\NN,<^\NN)$ is not well founded.
        \item The \textit{strong well ordering number} (essentially \cite{shelah2012nice}) of a logic $\LL$, if it exists, is the least ordinal $\alpha$ such that for any $\LL$-sentence $\phi$ in a signature containing $P$ and $<$ the following holds: if $\phi$ has a model $\MM$ such that
        \[
        (P^\MM,<^\MM)\cong(\alpha,\in),
        \]
        then $\phi$ also has a model $\NN$ such that
        \begin{itemize}
            \item $(P^\NN,<^\NN)$ is not well founded.
            \item There are $\{a_n\}_{n\in\omega}\subseteq P^\NN$ such that $P^\NN=\bigcup_{n\in\omega}\{b:b<^\NN a_n\}$.
        \end{itemize}
    \end{enumerate}
\end{de}

\begin{de}\spa 
\begin{enumerate}
    \item A logic or a fragment $\LL$ has the \textit{Löwenheim-Skolem-Tarski property} at a cardinal $\theta$ if for every model $\MM$ in a signature $\tau$ of size $\leq\theta$ there is $\NN\elem_\LL\MM$ of size $\leq\theta$.
    \item A logic or a fragment $\LL$ has the \textit{Tarski Union property} if whenever $(\MM_n)_{n\in\omega}$ is an $\LL$-elementary chain, then ${\MM_m\elem_\LL\bigcup_{n\in\omega}\MM_n}$ for each $m\in\omega$.
\end{enumerate}
    
\end{de}

\begin{rmk}
    Cartagena logic has Löwenheim-Skolem-Tarski property at $\kappa$ (Löwenheim-Skolem-Tarski Theorem \ref{LST}) and it satisfies the Tarski Union Property (Union Lemma \ref{union}). Furthermore, every fragment $\LL$ of size $<\kappa$ has Löwenheim-Skolem-Tarski property at some cardinal $<\kappa$, and satisfies Tarski Union Property (by the same theorems).
\end{rmk}

Following Definition \ref{fragelem}, we say that a map $\pi:\NN\to\MM$ is \textit{$\LL$-elementary} if
\[
\pi[\NN]\elem_\LL\MM.
\]

\begin{thm}[Strong undefinability of well order]\label{sudwo}
    Let $\theta$ be a cardinal and let $\LL$ be a fragment\footnote{Or an abstract logic.}. If $\LL$ satisfies
    \begin{itemize}
        \item the Löwenheim-Skolem-Tarski property at $\theta$, and
        \item the Tarski Union property,
    \end{itemize}
    then the \textit{strong well ordering number of $\LL$} is at most $(2^\theta)^+$.
\end{thm}
\begin{proof}
We show that for any model $\MM$ with $(P^\MM,<^\MM)=((2^\theta)^+,\in)$ there is $\NN\equiv_\LL\MM$ such that
\begin{enumerate}
    \item $(P^\NN,<^\NN)$ is not well founded.
    \item There are $\{a_n\}_{n\in\omega}\subseteq P^\NN$ such that $P^\NN=\bigcup_{n\in\omega}\{b:b<^\NN a_n\}$.
\end{enumerate}

Let $\MM$ be a model with $(P^\MM,<^\MM)=((2^\theta)^+,\in)$. We assume for simplicity that the signature has size $\leq\theta$. We proceed in four claims.

\begin{claim*}[\textbf{A}]  For each $n$, there is a collection of models $\{\NN^i_n:i\in(2^\theta)^+\}$  such that:
\begin{enumerate}
    \item\label{eka1} $i\in\NN^i_n\subseteq\NN^i_{n+1}\elem_\LL\MM$.
    \item\label{toka1} $|\NN^i_n|=\theta$.
    \item\label{kolmas1} $\sup\left(P^{\NN^i_n}\right)<\sup\left(P^{\NN^i_{n+1}}\right)$.
\end{enumerate}
\end{claim*}
This can be constructed by recursion on $n$, repeatedly using the Löwenheim-Skolem-Tarski property of $\LL$ at $\theta$.

\begin{claim*}[\textbf{B}]
    There are sets $X_0\supseteq X_1\supseteq X_2\supseteq\dots$ such that each $X_n\subseteq (2^\theta)^+$ is cofinal, and for all $i,j\in X_n$ there is an isomorphism
    \[
    \NN^i_n\cong\NN^j_n
    \]
    that maps $i\mapsto j$.
\end{claim*}

Again, these sets $X_n$ can be constructed by recursion on $n$, at each step relying on a counting argument: each structure $(\NN^i_n,i,\alpha^i_n)$ has size $\theta$ and there are only $2^\theta$ many pairwise non-isomorphic structures of size $\theta$.

\begin{claim*}[\textbf{C}]
    There are models
    \[
    (\NN_n,i_n,a_n)_{n\in\omega}
    \]
    and $\LL$-elementary maps
    \[
    \begin{tikzcd}
        &\NN_0\arrow{r}{\pi_0} &\NN_1\arrow{r}{\pi_1} &\NN_2\arrow{r}{\pi_2} &\dots
    \end{tikzcd}
    \]
    such that the elements $i_n,a_n\in P^{\NN_n}$ satisfy
    \begin{enumerate}
        \item $i_n\in X_n$.
        \item $i_{n+1}<^{\NN_{n+1}}\pi_n(i_n)$.
        \item $\sup\left(\pi_n[P^{\NN_n}]\right)<^{\NN_{n+1}}a_{n+1}$.
    \end{enumerate}
\end{claim*}
\begin{proof}[Proof of Claim (C)]
    
Each model $\NN_n$ is carefully chosen among the models 
\[
\{\NN^i_n:i\in X_n\},
\] for each $n$, again recursively:

For $n=0$, pick any $i\in X_0$ and let $\NN_0:=\NN_0^i$,  $i_0:=i$ and $a_0:=i$. 

For $n+1$, assume that $(\NN_n,i_n,a_n)$ has been defined and $i_n\in X_n$. We find $(\NN_{n+1},i_{n+1},a_{n+1})$ and the map $\pi_n:\NN_n\to\NN_{n+1}$.
Pick ordinals $i,j\in X_{n+1}$ such that $i<j$ and 
\[
i\in P^{\NN^j_{n+1}}\cap X_n.
\] 
This is possible because the sets $X_{n+1}$ are cofinal in $(2^\theta)^+$.
Now, since $j\in X_{n+1}\subseteq X_n$, there is an isomorphism $\NN_n\cong \NN^j_n$. Furthermore, $\NN^j_n\elem_{\LL}\NN^j_{n+1}$, because $\NN^j_n\subseteq\NN^j_{n+1}$ and both are $\LL$-elementary submodels of $\MM$. We let $\pi_n$ to be the composition of the isomorphism and the inclusion.
    \[
    \begin{tikzcd}
        &\NN^j_n\arrow{r}{\elem_\LL}&\NN^j_{n+1} \\
        &\NN_n\arrow{u}{\cong}\arrow[swap]{ur}{\pi_n} &
    \end{tikzcd}
    \]
We have:
\begin{itemize}
    \item $i\in X_{n+1}$.
    \item $\pi_n(i_n)=j>i$.
    \item $\sup\left(\pi_n[P^{\NN_n}]\right)=\sup\left( P^{\NN^j_n}\right)<\sup\left(P^{\NN^j_{n+1}} \right)$.
\end{itemize}
We let $\NN_{n+1}:=\NN^j_{n+1}$,  $i_{n+1}:=i$, and let $a_{n+1}\in P^{\NN^j_{n+1}}$ be any element such that $\sup\left(P^{\NN^j_n}\right)<a_{n+1}$. Then $\pi_n:\NN_n\to\NN_{n+1}$ and $(\NN_{n+1}, i_{n+1}, a_{n+1})$ are as wanted.
\end{proof}

Let $\NN$ be the direct limit
    \[
    \NN:=\dirlim\left(\NN_n,\pi_n:n<\omega\right).
    \]

\begin{claim*}[\textbf{D}]\spa 

    \begin{enumerate}
        \item $\NN\equiv_\LL\MM$.
        \item $(P^\NN,<^\NN)$ is not well founded.
        \item\label{third} There are $\{a_n\}_{n\in\omega}\subseteq P^\NN$ such that $P^\NN=\bigcup_{n\in\omega}\{b:b<^\NN a_n\}$.
    \end{enumerate}
\end{claim*}

Up to isomorphic correction, we may assume that each $\pi_n$ is in fact inclusion, $(\NN_n)_{n\in\omega}$ is an elementary chain and $\NN$ is its union.

By Tarski Union Property, $\NN_n\elem_\LL\NN$ for each $n$. This in particular implies $\NN\equiv_\LL\MM$. The fact that $(P^\NN,<^\NN)$ is not well founded is witnessed by the $i_n$'s:
\[
i_0>i_1>i_2>\dots.
\]
Finally, the item \ref{third} follows from the fact that
\[
P^{\NN_n}\subseteq\{b:b< a_n\},
\]
for each $n$.
This suffices to see that $\NN$ is as wanted.

\end{proof}

Applying Theorem \ref{sudwo} to each fragment $\LL_\phi$, we obtain:

\begin{cor}
    The (strong) well ordering number of Cartagena logic $\LL^c_\kappa$ is $\kappa$.
\end{cor}

\section{Expressive power of Cartagena logic}\label{sectexp}

The aim of this section is to compare Cartagena logic with the traditional infinitary logics $\LL_{\kappa\omega}$ and $\LL_{\kappa\kappa}$ and with Shelah's logic $\LL^1_\kappa$. In order to shed light on the expressive power of Cartagena logic, we will start by exhibiting five natural classes of models that are definable in it. Then, we will define \textit{$\omega$-covering property up to a cardinal $\theta$} and show that the class of models with this property is also definable, for each $\theta<\kappa$. This is crucial, since as a corollary, we obtain that the $\Delta$-closure of Cartagena logic is $\LL^1_\kappa$. The rest of the section is devoted to explaining this corollary.

We start with examples. Throughout this section, we assume that $\kappa$ is an uncountable cardinal satisfying $\kappa=\beth_\kappa$.

\subsection{Examples of expressive power}\label{examples}\spa 

\vv

\begin{thm}\label{expressive}
    The following classes of structures are definable in the Cartagena logic $\LL^c_\kappa$, for each cardinal $\theta<\kappa$:
\begin{enumerate}
    \item\label{card} Models of cardinality $\theta$, as well as models with a predicate or a definable subset of size $\theta$.
    \item Graphs with a clique of size $\theta$.
    \item Graphs of size $\theta$ that admit an $\omega$-coloring.
    \item For each cardinal $\theta<\kappa$ of uncountable cofinality: $\theta$-Aronszajn trees.
    \item Partially ordered models with an uncountable descending chain.
\end{enumerate}
\end{thm}
\begin{proof} We explicit a sentence defining each class.

    \begin{enumerate}
        \item \textbf{Cardinalities}

        Let $\theta<\kappa$ be a cardinal of uncountable cofinality.
        The sentence
        \[
        \forall\x_\theta\bigwedge_{f:\theta\to\omega}\bigdvee_{X\in\WW_f}\left(\bigdvee_{u\in[X]^\theta}\bigwedge_{\{i,j\}\in [u]^2}x_i=x_j\right)
        \]
        defines the class of structures of size $<\theta$, i.e. the class $\{\MM:|\MM|<\theta\}$. Furthermore, it is a Cartagena sentence, as the subformula \[
        \psi:=\bigdvee_{u\in[X]^\theta}\bigwedge_{\{i,j\}\in [u]^2}x_i=x_j
        \] 
        is upwards correct for $X$ and as for each $A\subseteq\theta$,
        \begin{align*}
            \FreeVar(\psi[A/X])&=\FreeVar\left(\bigvee_{u\in[A]^\theta}\bigwedge_{\{i,j\}\in [u]^2}x_i=x_j\right)\\
            &\subseteq\{x_i:i\in A\}.
        \end{align*}
        Similarly, given a formula $\phi(x)$, the sentence
        \[
        \forall\x_\theta\bigwedge_{f:\theta\to\omega}\bigdvee_{X\in\WW_f}\left(\bigdvee_{u\in[X]^\theta}\left(\bigwedge_{i\in u}\phi(x_i)\to\bigwedge_{\{i,j\}\in [u]^2}x_i=x_j\right)\right)
        \]
        defines the class of models in which the subset defined by $\phi(x)$ has size $<\theta$, i.e. the class $\{\MM:|\phi(\MM)|<\theta\}$.

        The class of models of cardinality at least $\theta$,  where $\theta$ has countable cofinality, as well as models of size at most $\theta$ or exactly $\theta$  can be defined by a Boolean combination of above sentences.

        \item \textbf{Graphs with a large clique} 

        The class of graphs with a clique of size $\theta$, for a  cardinal $\theta<\kappa$ of uncountable cofinality, is defined by the sentence
\[
\exists\x_\theta\bigwedge_{f:\theta\to\omega}\bigdvee_{X\in\WW_f}\left(\bigdvee_{u\in[X]^\theta}\bigwedge_{\{i,j\}\in[u]^2}E(x_i,x_j)\right)
\]
together with the conjunction of the graph axioms. For a cardinal of countable cofinality the defining sentence is a conjunction of the above ones.

        \item \textbf{Graphs with an $\omega$-coloring}

        The class of graphs of size $\theta$ that are colorable by countably many colors is defined by the sentence
        \[
            \forall\x_\theta\bigvee_{f:\theta\to\omega}\bigdwedge_{X\in\WW_f}\left(\bigdwedge_{\{i,j\}\in[X]^2}\neg E(x_i,x_j)\right)
        \]
        together with graph axioms and the sentence that defines models of size $\theta$.

        \item \textbf{$\theta$-Aronszajn trees} 
        
        A partial order $T=(T,<)$ is a \textit{$\theta$-Aronszajn tree} if every level $T_\alpha$ of $T$ has size $<\theta$, for every $t\in T$, the set $\{s:s<t\}$ is well ordered with order type $<\theta$ and there is no linearly ordered subset $b\subseteq T$ of size $\theta$.
        For $\theta<\kappa$, the class of $\theta$-Aronszajn trees is defined in Cartagena logic by the conjunction of axioms for partial order, the sentence that defines models of size $\theta$, and the following:
        \begin{enumerate}
            \item\label{a3} Each element has well ordered predecessors of order type $<\theta$: 
            \[
            \forall x\bigvee_{\alpha<\theta} \q \textit{The set of predecessors of }x\textit{ has order type }\alpha\q.
            \]
            \item\label{a4} Levels of size $<\theta$:
            \[
            \bigwedge_{\alpha<\theta}\q\textit{The }\alpha\textit{th level }T_\alpha\textit{ has size }<\theta\q.
            \]
            \item\label{a5} No long branch:
            \[
            \forall\x_\theta\bigwedge_{f:\theta\to\omega}\bigdvee_{X\in\WW_f}\bigdvee_{\{i,j\}\in[X]^2}(x_i\bot x_j\vee x_i=x_j).
            \]
        \end{enumerate}
        The sentence \ref{a3} is an $\LL_{\kappa\omega}$ sentence. For details, see IX.1.2 from \cite{Barwise1985-BARML-8}. For the same reason, each $\alpha$th level $T_\alpha$ is an $\LL_{\kappa\omega}$-definable hence $\LL^c_\kappa$-definable set, thus application of item (\ref{card}) gives that \ref{a4} is a Cartagena sentence. The sentence \ref{a5} is a Cartagena sentence since the subformula $\psi:=\bigdvee_{\{i,j\}\in[X]^2}(x_i\bot x_j\vee x_i=x_j)$ is upwards correct with respect to $X$ and for each $A\subseteq\theta$,
        \[
        \FreeVar(\psi[A/X])\cap\{x_i:i\in\theta\}\subseteq\{x_i:i\in A\}.
        \]

        \item \textbf{Uncountable descending chain}
        
        Let $\theta:=\omega_1$. The sentence
        \[
        \exists\x_\theta\bigwedge_{f:\theta\to\omega}\bigdvee_{X\in\WW_f}\bigdvee_{u\in[X]^\theta}\bigwedge_{\substack{\{i,j\}\in [u]^2,\\
        i<j}}x_j<x_i
        \]
        defines the class of models with an uncountable descending chain.

    \end{enumerate}
\end{proof}

Next we describe a covering property that will be used when examining the distance of Cartagena logic from Shelah's logic $\LL^1_\kappa$.

\begin{de}\label{covering}
    A model $\MM$ in a signature containing a binary predicate symbol $E$ has the \emph{$\omega$-covering property up to $\theta$ (with respect to $E$)} if
    \begin{enumerate}
        \item For every subset $A\subseteq\MM$ of size $\leq\theta$ there are $\{a_n\}_{n\in\omega}\subseteq\MM$ such that the set of $E$-predecessors of the $a_n$'s cover $A$:
        \[
        A\subseteq\bigcup_{n\in\omega}\{b:E^\MM(b,a_n)\}.
        \]
        \item For each $b\in\MM$,
        \[
        |\{c\in\MM:E^\MM(c,b)\}|\leq\theta.
        \]
    \end{enumerate}
\end{de}

Heuristically, the $\omega$-covering property states that every small subset can be covered with countably many designated small sets.

\begin{ex} The real line $(\R,<)$ has the $\omega$-covering property up to $|\R|$.
\end{ex}

\begin{ex} Every ordinal of countable cofinality, as a structure $(\alpha,\in)$, has countable covering property up to $|\alpha|$.
\end{ex}

\begin{prop}\label{coveringdef}
    The class of models with the $\omega$-covering property up to $\theta$ is definable in Cartagena logic, for each $\theta<\kappa$.
\end{prop}
\begin{proof} Write
\begin{align*}
    \Theta_\theta:= \forall\x_\theta\bigvee_{f:\theta\to\omega}\bigdwedge_{X\in\WW_f}\exists y\left(\bigdwedge_{i\in X}E(x_i,y)\right).
\end{align*}

    For each $\theta<\kappa$, the sentence $\Theta_\theta$ is a Cartagena sentence, as the subformula
    \[
    \psi(\x_\theta,X):=\exists y\left(\bigdwedge_{i\in X}E(x_i,y)\right)
    \]
    is downwards correct for $X$ and satisfies that for each $A\subseteq\theta$
    \[
    \FreeVar\left(\psi[A/X]\right)\cap\{x_i:i\in\theta\}\subseteq\{x_i:i\in A\}.
    \]
    Furthermore, $\Theta_\theta$ together with the sentence
    \[
    \forall x\spa \textit{``the set of }E\textit{-predecessors of }x\textit{ has size }\leq\theta",
    \]
    which is a Cartagena sentence by \ref{expressive} (\ref{card}),
    defines the class of models with $\omega$-covering property up to $\theta$.
\end{proof}

We will now describe Shelah's logic $\LL^1_\kappa$, to be able to compare Cartagena logic with it.

\subsection{Shelah's logic $\LL^1_\kappa$}\label{sectshlogic}\spa 

\vv 

This section is devoted to a discussion of the logic $\LL^1_\kappa$. We first give the relevant game, \textit{delayed game} (called as such in \cite{velickovicVaananen}), and then derive the logic $\LL^1_\kappa$ out from it.

The delayed game is a more advanced version of the Cartagena game. Similarly to Cartagena game, player $\1$ picks sets of size $<\kappa$ and a descending sequence of ordinals, and player $\2$ has to map parts of the sets in the opposite model. As in Cartagena game, she partitions the set picked by player $\1$ with $\omega$.  However, opposed to Cartagena game in which player $\1$ chooses the piece she has to map, in delayed game she has to map the first piece of the partition. And on the next round, she has to map the second piece of the same partition, etc. At each round, player $\1$ introduces a new large set, so player $\2$ has an ever-increasing  amount of pieces to map. As the game clock renders the game finite, she will in fact never be forced to map all of any of the large sets - only a finite initial segment of the partition of each. Our presentation here has a height function $\heit$, which is used to indicate the pieces player $\2$ must map at each state: she maps those that have height $0$ at that state.

Inspired from the delayed game, originally Cartagena game was defined slightly differently from its current definition: first player $\1$ played a set, then player $\2$ played a partition of the set, then player $\1$ played a piece of the partition, and finally player $\2$ had to play a partial isomoprhism of this single piece. 

We changed large sets into long tuples, which allowed us to change the order of moves as they are in the present paper: player $\2$ simultaneously plays her partition and a long tuple in the opposite model. This change of order was a key for finding the syntax for Cartagena logic. Analogous move does not seem possible in the delayed game, a fact that heavily complicates the search for syntax for $\LL^1_\kappa$. 

We now define the delayed game.

\vv

In the following, we denote $\fld(\pi)=\dom(\pi)\cup\ran(\pi)$, and $\dot -$ denotes truncated subtraction.

\begin{de}[Delayed game $\DG^\beta_\lambda$]\label{shgame} Let $\lambda$ be a cardinal, let $\MM$ and $\NN$ be structures in a same signature with pairwise disjoint domains, and let $\beta$ be an ordinal. We define the \textit{delayed game} of height $\beta$ 
\[
\DG_\lambda^\beta(\MM,\NN).
\]

The \textit{states} of this game are triples $(\alpha,\pi,\heit)$, where $\alpha\leq\beta$ is an ordinal, ${\pi:\MM\to\NN}$ is a partial isomorphism and $\heit:\MM\cup\NN\to\omega$ is a partial function such that ${\fld(\pi)\subseteq\dom(\heit)}$.

\textbf{Starting state:} The starting state is $(\beta,\emptyset,\emptyset)$.

\textbf{Further states:} At state $(\alpha,\pi,\heit)$:
\begin{enumerate}
    \item Player $\1$ picks an ordinal $\alpha'<\alpha$ and a set $A\in[\MM]^{\leq \theta}\cup[\NN]^{\leq\theta}$, for some cardinal $\theta<\lambda$.
    \item Player $\2$ picks partial functions
    \begin{align*}
        & \pi':\MM\to\NN\\
        & \heit':\MM\cup\NN\to\omega
    \end{align*}
    such that:
    \begin{itemize}
        \item $\pi'$ is a partial isomorphism extending $\pi$.
        \item If $\heit'(a)=0$, then $a\in\fld(\pi')$.
        \item $A,\dom(\heit)\subseteq\dom(\heit')$.
        \item For all $a\in\dom(\heit)$, $\heit'(a)=\heit(a)\dot -1$.
    \end{itemize}
\end{enumerate}

The next state is $(\alpha', \pi',\heit')$.
    
\noindent The player who first cannot move loses.
\end{de}

\begin{de}
    We define $\sim^\beta_\lambda$ to be the transitive closure of the relation
    \[
    \textit{Player }\2\textit{ has a winning strategy in }\DG_\lambda^\beta(\MM,\NN).
    \]
\end{de}
It is not yet known whether the game itself is already transitive, so we have to content ourselves to taking the transitive closure.
We are now ready to define Shelah's logic $\LL^1_\kappa$.

\begin{de}[Logic $\LL^1_\kappa$] Let $\kappa$ be an uncountable cardinal such that $\kappa=\beth_\kappa$ and let $\tau$ be a signature of size $<\kappa$.\footnote{For simplicity. The general case can be found in \cite{shelah2012nice}.}
\begin{enumerate}
    \item A \textit{sentence} $\phi\in\LL^1_\kappa(\tau)$ is a class of $\tau$-structures which is closed under the relation $\sim^\beta_\theta$, for some $\beta,\theta<\kappa$. 
    \item For an expansion of a $\tau$-structure $\MM$ and a sentence $\phi\in\LL^1_\kappa(\tau)$:
    \[
    \MM\models\phi:\iff\MM\rest\tau\in\phi.
    \]
\end{enumerate}
\end{de}

We state a Lindstr\" om-style characterization without proof:

\begin{thm}[First Characterization Theorem for $\LL^1_\kappa$, Shelah, \cite{shelah2012nice}]\label{firstcha} For an uncountable cardinal $\kappa$ such that $\kappa=\beth_\kappa$, the logic $\LL^1_\kappa$ has strong well ordering number $\kappa$, Löwenheim-Skolem number\footnote{The Löwenheim-Skolem number of a logic, if exists, is the least cardinal $\lambda$ such that every sentence that has a model has a model of size $<\lambda$.} $\kappa$ and is a maximal such logic above $\LL_{\kappa\omega}$.
\end{thm}

Next we discuss the $\Delta$-closure of a logic, which is a weak form of interpolation, and observe that $\Delta(\LL^c_\kappa)=\LL^1_\kappa$.

\subsection{$\Delta$-closure of $\LL^c_\kappa$}\spa 

\vv

Logics are compared with respect to their expressive power. We write \[
\LL\leq\LL'
\] 
if for every $\LL$-sentence $\phi$ there is an $\LL'$-sentence $\psi$ with $\Mod(\phi)=\Mod(\psi)$.

The next easy lemma shows that $\LL^1_\kappa$ is at least as strong as Cartagena logic $\LL^c_\kappa$.

\begin{lem}\label{comparison}
    If player $\2$ has a winning strategy in $\DG_\kappa^\beta(\MM,\NN)$,
    then she has a winnning strategy in $\G_\kappa^{\beta\cdot\omega}(\MM,\NN)$.
    In particular:
    \[
    \LL^c_\kappa\leq\LL^1_\kappa.
    \]
\end{lem}

In fact, it will turn out in a moment that $\LL^1_\kappa$ is strictly stronger than Cartagena logic. We first discuss the $\Delta$-closure operation.

\begin{de} Let $\KK$ be a class of models in a fixed signature $\tau$, which may be many-sorted, and let $\LL$ be a logic.
    \begin{enumerate}
        \item $\KK$ is \textit{definable} in $\LL$ if there is an $\LL$-sentence $\phi$ in the signature $\tau$ such that
        \[
        \KK=\Mod(\phi).
        \]
        \item $\KK$ is \textit{projective} in $\LL$ if there is an $\LL$-sentence $\phi$ in an expanded signature $\tau'\supseteq\tau$, which may have new sorts, such that
        \[
        \KK=\{\MM\rest\tau:\MM\in\Mod(\phi)\}.
        \]
        \item $\LL$ is \textit{$\Delta$-closed} if every class of models which is both projective and co-projective in $\LL$ is in fact definable in $\LL$.
    \end{enumerate}
\end{de}

The following is well defined (see \cite{Barwise1985-BARML-8}, II.7.2, especially Definition 7.2.3.):

\begin{de}
    The \textit{$\Delta$-closure} of a logic $\LL$ is the smallest logic $\Delta(\LL)\geq\LL$ which is $\Delta$-closed.
\end{de}

\begin{rmk}
    The $\Delta$-closure is a well defined closure operation
    \[
    \LL\mapsto\Delta(\LL)
    \]
    defined on every logic $\LL$. It preserves many model-theoretic properties, such as Löwenheim-Skolem number and (strong) well ordering number. It is not known to preserve the existence of generative syntax. See more in \cite{Barwise1985-BARML-8}, II.7.2.

\end{rmk}

\begin{rmk}\label{inter} A logic $\LL$ has \textit{interpolation} if for any sentences $\phi\in\LL(\tau)$ and ${\phi'\in\LL(\tau')}$, if $\phi\models\phi'$, then there is an interpolant $\psi\in\LL(\tau\cap\tau')$ with $\phi\models\psi\models\phi'$. Interpolation is equivalent to the following: \textit{any two disjoint projective model classes are separable by a definable model class}, which clearly is a strengthening of being $\Delta$-closed.
    
\end{rmk}

Again, we quote a theorem without a proof:

\begin{thm}[Shelah, \cite{shelah2012nice}] The logic $\LL^1_\kappa$ has interpolation.
\end{thm}

In particular, the logic $\LL^1_\kappa$ is $\Delta$-closed, by Remark \ref{inter}. Furthermore, it has a characterization in terms of $\Delta$-closure:

\begin{thm}[Second Characterization Theorem for $\LL^1_\kappa$, Shelah, \cite{shelah2012nice}]\label{secondcha} The logic $\LL^1_\kappa$ is the minimal logic above $\LL_{\kappa\omega}$ which is $\Delta$-closed and in which the class of models with $\omega$-covering property up to $\theta$ is definable, for each $\theta<\kappa$.
\end{thm}

The class of models with $\omega$-covering property up to $\theta$ is definable in Cartagena logic, for each $\theta<\kappa$, by Proposition \ref{coveringdef}. Thus application of Theorem \ref{secondcha} gives:

\begin{cor}
    $\Delta(\LL^c_\kappa)=\LL^1_\kappa$.
\end{cor}

We will now find two structures that are elementary equivalent in Cartagena logic but not in $\LL^1_\kappa$. This will show that $\LL^c_\kappa<\LL^1_\kappa$. These two structures are the real line $(\R,<)$ and the real line without zero $(\R-\{0\},<)$. In \cite{velickovicVaananen}, the following is proved:

\begin{lem}[Väänänen, Veličkovi\' c, \cite{velickovicVaananen}]\label{vel} $(\R,<)\not\equiv_{\LL^1_\kappa}(\R-\{0\},<)$.
    
\end{lem}

There remains to show the following:

\begin{lem} $(\R,<)\equiv_{\LL^c_\kappa}(\R-\{0\},<)$.
\end{lem}
\begin{proof} We only sketch the proof. It suffices to describe a winning strategy for player $\2$ in the game
        \[
        \G_\kappa^\beta((\R,<),(\R-\{0\},<)),
        \]
        where $\beta$ is an arbitrary ordinal. The argument is to show by induction that player $\2$ can maintain the following condition at each state $(\alpha,\pi)$:
        \begin{itemize}
            \item[($*$)] \textit{There are finitely many pairwise disjoint closed intervals $I_0,\dots,I_n\subseteq\R$ and $J_0,\dots,J_n \subseteq\R-\{0\}$ such that
            \begin{enumerate}
                \item $\dom(\pi)=\bigcup_{k\leq n} I_k$.
                \item For each $k\leq n$,
                    \[
                        \pi\rest I_k:(I_k,<)\cong (J_k,<).
                    \]
            \end{enumerate}}
        \end{itemize}
    
\end{proof}

To summarize: 

\begin{cor}\spa 
    \begin{enumerate}
        \item $\LL^c_\kappa<\LL^1_\kappa$.
        \item $\Delta(\LL^c_\kappa)=\LL^1_\kappa$.
        \item Cartagena logic is not $\Delta$-closed and does not have interpolation.
        \item The $\Delta$-closure of Cartagena logic has interpolation.
    \end{enumerate}
\end{cor}

\section*{Further work}

The goal of our work was to find a generative syntax for the logic $\LL^1_\kappa$, and therefore understand better the general question of when and how is it possible to derive a syntax from a game. Partial results were achieved. The method used in the present paper, as it is, does not directly give a syntax for $\LL^1_\kappa$, because the delayed game (\ref{shgame}) is more involved. These subtleties were discussed in detail in the beginning of Section \ref{sectshlogic}. 

However, we believe that the strategy presented here can be pushed further to build a simple generative syntax for $\LL^1_\kappa$ too, which will elucidate further the general problem of deriving syntax from a game. We also have hopes that the Boolean extension $\LL^{\Bool}_{\kappa\kappa}$ could prove useful elsewhere.

\bibliographystyle{plain}
\bibliography{bib}

\end{document}